\documentclass [11pt,reqno]{article}

\usepackage{enumerate,amsthm,amsmath,amsfonts,amssymb,pstricks,hyperref}
 \numberwithin{equation}{section}

\makeatletter \xdef\qedbuit{\qed}

\newcommand{\TeoremaAmbFinalMarcat}[1]{%
  \expandafter\gdef\csname end#1\endcsname{\qedbuit\@endtheorem}}
\newtheorem{theo}{Theorem}[section]

\newtheorem{cory}[theo]{Corollary}

\theoremstyle{definition}
 \TeoremaAmbFinalMarcat{defi}
 \TeoremaAmbFinalMarcat{ex}
\newtheorem{rem}[theo]{Remark} \TeoremaAmbFinalMarcat{rem}
 \TeoremaAmbFinalMarcat{nota}
 \TeoremaAmbFinalMarcat{prob}
 \TeoremaAmbFinalMarcat{ques}
\newenvironment{proclama}[1]{\par\vspace{\topsep}{\bf  #1}
 \begin{em}}{\end{em}\par\vspace{\topsep}}

\newcommand{\start}[2]{\begin{#1}\label{#2}}

\newcommand{\theoc}[1]{Theorem~\ref{#1}}
\newcommand{\propc}[1]{Proposition~\ref{#1}}
\newcommand{\coryc}[1]{Corollary~\ref{#1}}
\newcommand{\lemc}[1]{Lemma~\ref{#1}}

\newcommand{\remc}[1]{Remark~\ref{#1}}
\newcommand{\exc}[1]{Example~\ref{#1}}
\newcommand{\figc}[1]{Figure~\ref{#1}}

\newcommand{\refc}[1]{~\ref{#1}}
\newcommand{\notc}[1]{Notation~\ref{#1}}

\def\@enum@{\list{\csname label\@enumctr\endcsname}%
           {\usecounter{\@enumctr}\def\makelabel##1{\hss\llap{##1}}
           \itemsep=2pt\parsep=0pt\topsep=3pt plus 1pt minus 1 pt}}

\newenvironment{romlist}{\enumerate[(i)]}{\endenumerate}
\newenvironment{numlist}{\enumerate[(1)]}{\endenumerate}

\def\unitC{1}

\def\map#1#2#3{\mbox{${#1}\colon {#2} \longrightarrow {#3}$}}

\newcommand{\SI}{\ensuremath{\mathbb{S}}^1}
\newcommand{\Z}{\ensuremath{\mathbb{Z}}}
\newcommand{\T}{\ensuremath{\mathbb{T}}}
\newcommand{\R}{\ensuremath{\mathbb{R}}}

\begin{document}

\title{Minimal intersection of curves on surfaces}\date{}
\author{Moira Chas
\thanks{This work is partially supported in part by NSF grant 1034525-1-29777.}
}


\maketitle

\begin{abstract} This paper is a consequence of the close connection between combinatorial group theory and the topology of surfaces.

In the eighties Goldman discovered a Lie algebra structure on the vector space generated by the free homotopy classes of oriented curves on an oriented surface.  The Lie bracket [a,b] is defined as the signed sum over the intersection points  of  a and b of the loop product of at the intersection points.

If one of the classes has a simple representative  we give a combinatorial group theory  description of the terms of  the Lie bracket and  prove that this bracket has as many terms, counted with multiplicity, as the minimal  number of intersection points of a and b. In other words the bracket with a simple element has no cancellation and determines minimal intersection 
numbers. We show that analogous results hold for the Lie bracket (also discovered by Goldman) of unoriented curves. We give three applications: a factorization of Thurston's map defining the boundary of Teichm\"{u}ller space,  various decompositions of the underlying vector space of conjugacy classes into ad invariant subspaces and a connection between bijections of the set of conjugacy classes of curves  on a surface preserving the Goldman bracket and the mapping class group.

\medskip

{\footnotesize
\noindent \emph{{\normalfont 2000}\,Mathematics Subject Classification.} Primary: 57M99,  20E06.

}

\end{abstract}

\section{Introduction}

Let $a$ and $b$ denote isotopy classes of embedded closed curves on a surface $\Sigma$. Denote
by $i(a, b)$ the minimal possible number of intersection points 
 of curves representing $a$ and $b$, where the intersections are counted with multiplicity. The
function $i(a, b)$ plays a central role in Thurston's work on low dimensional topology (see, for
instance, \cite{T}, \cite{flp} and \cite{hp}.) Let $[a, b]$ denote the Lie bracket on the vector
space of the free homotopy classes of all essential directed closed curves on $\Sigma$. This Lie
bracket  originated from Wolpert's  cosine formula, Thurston's  earthquakes in Teichm\"{u}ller space and
Goldman's study of Poisson brackets (see \cite{g}.) The Lie algebra of oriented curves has been generalized using 
the loop product to a string bracket on 
the reduced $\SI$-equivariant homology of the free loop space of any oriented manifold. 
(see \cite{CS1} and \cite{CS2}.)

The main  purpose of this paper is to prove that both, the Goldman Lie bracket on oriented curves and the Goldman Lie bracket on unoriented curves  "count" the 
minimal number of intersection points of two simple curves.

In order to give a more precise statement of our results, let us review the definition of the Lie algebra for oriented curves discovered by Goldman: given two such free 
homotopy classes of
directed curves on an orientable surface, take two representatives intersecting each other only in
transverse double points. Each one of the intersection points will contribute a term to the geometric formula for the 
bracket. Each of these terms is defined as follows:  take the conjugacy class of the curve obtained
by multiplying the two curves at the intersection point and adjoin a negative sign if the orientation
given by the ordered tangents at that point is different from the orientation of the surface. The bracket on $W,$ the 
vector space of free homotopy classes is the bilinear extension of this construction. It is remarkable that this 
construction is well defined and satisfies skew-symmetry and   Jacobi .

The bracket of unoriented curves can be defined on the subspace $V$ of $W$ fixed by the operation of reversing 
direction, as the restriction of the bracket. The subspace $V$ is generated by elements of $W $ of the form   $a+
\overline{a}$ where $a$ is a basis element and where $\overline{a}$ denotes $a$ with opposite direction.  Let us identify unoriented curves up to 
homotopy with these expressions  $a+\overline{a}$.
The Lie bracket of two unoriented curves, $a+\overline{a}$ and $b+\overline{b}$, is then defined geometrically by
$([a,b]+\overline{[a,b]})+([a,\overline{b}]+\overline{[a,\overline{b}]}),$ which equals $[a+\overline{a},b+\overline{b}]$ 
using $\overline{[a,b]}=[\overline{a},\overline{b}]$. An \emph{unoriented term}  is a term  of the form $c(z +\overline
{z})$, where $c$ is an integer and $z$ is a conjugacy class, that is, an element of the basis of $V$ multiplied by an 
integer coefficient.

Since this Lie bracket uses the intersection points of curves, a natural problem was to study how
well it reflects the intersection structure. In this regard, Goldman \cite{g} proved the following
result.

\begin{proclama}{Goldman's Theorem} If the bracket of two free homotopy classes of curves (oriented or 
unoriented) is zero, and one of them has a simple representative, then the two classes have disjoint 
representatives.
\end{proclama}

Goldman's proof uses the Kerckhoff earthquake convexity property of Teichm\"{u}ller space \cite{Ke} and in
\cite{g} he wondered whether this topological result had a topological proof. In \cite{chas} we  gave
such a proof when $a$ was  a non-separating simple closed curve on a surface with non-empty
boundary.  Here we will give  combinatorial proof of our Main Theorem (see below), which  a generalization of  Goldman's result.

Now, suppose that $a$ can be represented by a simple closed curve
$\alpha$. Then for any free homotopy class $b$ there exists a representative that can be written as a
certain product which involves a sequence of elements in the
fundamental group or groups of the connected components of $\Sigma
\setminus \alpha$ (see Sections \ref{Amalgamated products} and
\ref{hnn extensions} for precise definitions.) The number of terms
of the sequence for $b$ with respect to $a$ is denoted by $t(a,b)$.

\begin{proclama}{Main Theorem} Let $a$ and $b$ be  two free homotopy
classes of directed curves on an orientable surface. If $a$ can be represented by a simple closed curve then the
following nonnegative integers are equal:
\begin{romlist}
\item The number of terms, counted with multiplicity, of the Goldman Lie bracket for oriented curves $[a,b]$ .
\item The number  divided by two of unoriented terms (of the form $x +\overline{x}$), counted with multiplicity, of 
the Goldman Lie bracket for unoriented curves $[a+\overline{a},b+\overline{b}]$.
\item The minimal number of intersection points of $a$ and $b$, $i(a,b)$.
\item The number of terms of the sequence for $b$ with respect to $a$,
$t(a,b).$
\end{romlist}
\end{proclama}

In particular, there is no cancellation of terms in the bracket of two curves, provided that one of them is simple.

As a consequence of our Main Theorem, we obtain the following result.

\begin{proclama}{Corollary of the Main Theorem} If $x$ and  $y$ are conjugacy classes of curves that can be represented by simple closed curves then 
$t(x,y)=t(y,x)$.
\end{proclama}

We obtain the following global characterization of conjugacy classes containing simple representatives. 

\begin{proclama}{Corollary of the Main Theorem} Let $a$ denote a free homotopy class of curves on a surface. Then the following statements are equivalent. \begin{numlist}
\item The class $a$ has a representative that is a power of a simple curve.
\item  For every free homotopy class $b$ the (oriented) bracket $[a,b]$ has many as  oriented terms   counted with multiplicity as the minimal intersection number of $a$ and $b$.
\item  For every free homotopy class $b$ the (unoriented) bracket $[a,b]$ has many as unoriented terms   counted with multiplicity as  twice the minimal intersection number of $a$ and $b$.
\end{numlist}
\end{proclama}

In \cite{ck}  a  local characterization of simple closed curves in terms of the Lie bracket will be given. Actually, the 
problem of characterizing algebraically embedded conjugacy classes was the original motivation of \cite{chas}, \cite
{CS1} and \cite{CS2}.

Here is a brief outline of the arguments we follow to prove the main theorem.

\begin{numlist}
\item The key
point is that when a curve is simple, we can apply either HNN
extensions or amalgamated products to  write elements of the
fundamental group of the surface as a product that involves certain sequences of
elements of subgroups which are the fundamental group of the connected components of the surface minus the simple curve. 
\item  Using combinatorial
group theory tools we show that if certain equations do not hold in an HNN extension or an
amalgamated free product then certain products of the sequences in (1) cannot be conjugate, (Sections \ref
{Amalgamated products} and \ref{hnn extensions}.)
\item We show that each of the terms of the bracket can be obtained by
inserting the conjugacy class of the simple curve in different
places of the sequences in (1) made circular. (Sections \ref{separating oriented} and 
\ref{non-separating case}.)
\item We show that the equations of (2) do not hold in the HNN extensions and amalgamated free products 
determined by a simple closed curve. (Section \ref{aux})
\end{numlist}

The Goldman bracket extends to higher dimensional manifolds as one of the String Topology operations.  Abbaspour \cite{hossein} characterizes hyperbolic three manifolds among closed three manifolds using the loop product which is another String Topology operation. Some of his arguments rely on the decomposition of the fundamental group of a manifold into amalgamated products based on torus submanifolds, and the use of this decomposition to give expressions for certain elements of the fundamental group, which in term, gives a way of computing the loop product.

\begin{figure}[htbp]\label{figure separating}
\begin{center}
\begin{pspicture}(12,5)
 \psccurve(1,2)(1,0.2)(6,1)(9,0.3)(11,3)(6,3)(1,4.4)
\pscurve(2,1.5)(2.7,1.2)(3.5,1.5)
\pscurve(2.2,1.4)(2.7,1.5)(3.2,1.34)

\rput(-0.5,1.58){\pscurve(2,1.5)(2.7,1.2)(3.5,1.5)
\pscurve(2.2,1.4)(2.7,1.5)(3.2,1.34) }
\rput(6,0.51){\pscurve(2,1.5)(2.7,1.2)(3.5,1.5)
\pscurve(2.2,1.4)(2.7,1.5)(3.2,1.34) }

\pscurve[linecolor=gray,linewidth=1pt](6,1)(6.15,2)(6,3)
\pscurve[linecolor=gray,linewidth=1pt,linestyle=dashed](6,1)(5.9,2)(6,3)

\psccurve[showpoints=false,linecolor=lightgray](1,3)(3,2.5)(9.2,1.5)(9.6,3)(6.2,2)(2,4)(2,1)(6.1,1.6)(10,2)(6.1,2.6)

\rput(6,3.3){$\chi$}
\end{pspicture}
\end{center}
\caption{A separating curve $\chi$ intersecting another curve}
\end{figure}
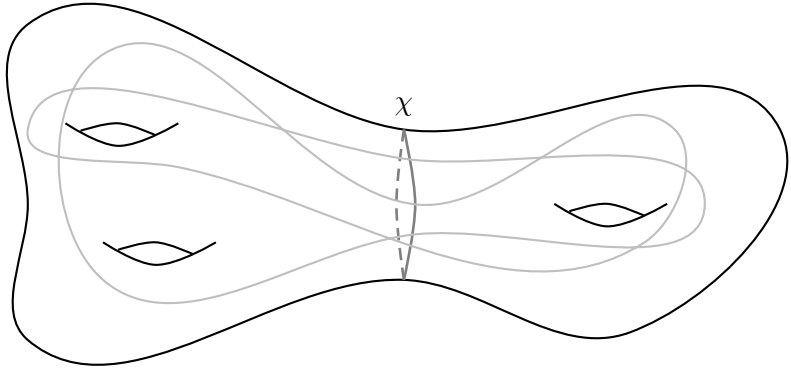

\begin{figure}[htbp]
\begin{center}
\begin{pspicture}(14,5)
 \psccurve(1,2)(1,0.2)(6,1)(9,0.3)(11,3)(6,3)(1,4.4)
\pscurve(2,1.5)(2.7,1.2)(3.5,1.5)
\pscurve(2.2,1.4)(2.7,1.5)(3.2,1.34)

\rput(-0.5,1.58){\pscurve(2,1.5)(2.7,1.2)(3.5,1.5)
\pscurve(2.2,1.4)(2.7,1.5)(3.2,1.34) }
\rput(6,0.51){\pscurve(2,1.5)(2.7,1.2)(3.5,1.5)
\pscurve(2.2,1.4)(2.7,1.5)(3.2,1.34) }

\pscurve[linestyle=dashed, linecolor=gray](8.8,2)(8.9,2.7)(8.8,3.4)
\pscurve[linecolor=gray](8.8,2)(8.7,2.7)(8.8,3.4)

\psccurve[showpoints=false,linecolor=lightgray](1,3)(3,2.5)(9.2,1.5)(9.6,3)(6.2,2)(2,4)(2,1)(6.1,1.6)(10,2)(6.1,2.6)

\rput(8.8,3.764){$\lambda$}

\end{pspicture} \caption{The intersection of a non-separating curve $\lambda$
and another curve}\label{figure non-separating}
\end{center}
\end{figure}
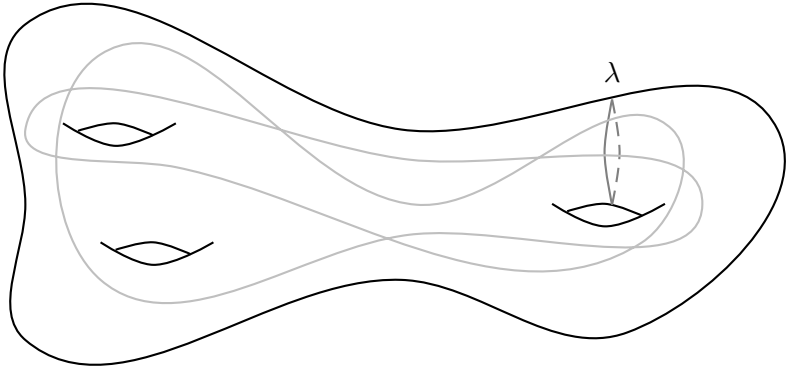

Here is the organization of this work: in Section \ref{Amalgamated products} we list the known results about amalgamated products of groups we will use, and we prove that certain equations do not hold in such groups. In Section \ref{separating oriented}, we apply the
results of the previous section to find a combinatorial description
of the Goldman bracket of  two oriented curves, one of them simple
and separating (see Figure \ref{figure separating}.) In Section
\ref{hnn extensions} we list results concerning HNN extensions and
we prove that certain equations do not hold in such groups. In
Section \ref{non-separating case}, we apply the results of Section
\ref{hnn extensions} to describe combinatorially the bracket of an
oriented non-separating simple  closed curve with another oriented
curve (see Figure \ref{figure non-separating}.) In Section \ref{aux}
we prove propositions about the fundamental group of the surface
which will be used to show that our sequences satisfy the hypothesis
of the main theorems of \ref{Amalgamated products} and \ref{hnn
extensions}. In Section \ref{conclusions} we put together most of the
above results to show that there is no cancellation in the Goldman
bracket of two directed curves, provided that one of them is simple.
In Section \ref{unoriented loops} we define the bracket of
unoriented curves and  prove that there is no cancellation if one of
the curves is simple. In Section \ref{examples} we exhibit some
examples that show that the hypothesis of one of the curves being
simple cannot be dropped (see Figure \ref{example}.) In the next three  sections, we exhibit some applications of our main
results. More precisely, in
Section \ref{applications 1}, we show how our main theorem yields a factorization of Thurston's map on the set of all simple conjugacy classes of curves on a surface, through the power of the vector space of all conjugacy classes to the set of simple conjugacy classes.
In Section~\ref{applications 2} we exhibit several partitions of the vector space generated by all conjugacy classes, which are invariant under certain Lie algebra operations. We conclude by showing in Section  \ref{applications 3}  that if a function on the set of conjugacy classes preserves the bracket, then is determined by an element of the mapping class group of the surface.
We conclude by stating some problems and open questions in Section \ref{questions}.

\textit{Acknowledgments:} This work benefitted from stimulating exchanges with 
Pavel Etingof and Dennis Sullivan and especially Warren Dicks who suggested an improved demonstration of \theoc{amalgamating with cyclic}.

\section{Amalgamated products}\label{Amalgamated products}

This section deals exclusively with results concerning Combinatorial
Group Theory. We start by stating  definitions and known results about amalgamating free products.  Using these 
tools, we prove the  main results of this section, namely,  Theorems \ref{amalgamating with cyclic} and \ref
{amalgamating inverse}. These two theorems state that certain pairs of elements of an amalgamating free product 
are not conjugate. 
By \theoc{bracket form}  if $b$ is an arbitrary conjugacy class and $a$ is a conjugacy class containing a separating 
simple representative, then the Goldman Lie bracket $[a, b]$ can be written as an algebraic sum $t_1 + t_2 + 
\cdots t_n$, with the following property: If there exist two terms terms $t_{i}$ and $t_{j}$ that cancel, then  the 
conjugacy classes associated to the terms $t_i$ and $t_j$  both satisfy the hypothesis of \theoc{amalgamating with 
cyclic} or both satisfy the hypothesis of \theoc{amalgamating inverse}.
This will show that the terms of the Goldman Lie bracket exhibited in \theoc{bracket form} are all distinct.

Alternatively, one could make use of the theory of groups acting on graphs (see, for instance,
 \cite{DD}) to prove Theorems \ref{amalgamating with cyclic} and \ref{amalgamating inverse}.

Let $C, G$ and $H$ be groups and let $\map{\varphi}{C}{G}$ and
$\map{\psi}{C}{H}$ be monomorphisms. We denote the \emph{free
product of $G$ and $H$ amalgamating the subgroup $C$ (and morphisms
$\varphi$ and $\psi$)} by  $G \ast_C H$. This group is defined as the quotient of the free product $G \ast H$ by the 
normal subgroup generated by $\varphi(c)\psi(c)^{-1}$ for all $c \in C$.
 (See \cite{h},\cite{LS}, \cite{MKS} or \cite{cohen} for detailed definitions.)

Since there are canonical injective maps from $C, G$ and $H$ to $G
\ast_C H$, in order to make the notation lighter we will work as if
$C$, $G$ and $H$ were included in $G\ast_C H$ .

\start{defi}{reduced form} Let $n$ be a non-negative integer. A finite sequence  $(w_1, w_2, \ldots,w_n)$ of 
elements  of $G \ast_C H$  is \emph{reduced} if the following
conditions hold,
\begin{numlist}
\item Each $w_i$ is in one of the factors, $G$ or $H$.
\item For each $i \in \{1,2,\ldots, n-1\}$, $w_i$ and $w_{i+1}$ are
not in the same factor.
\item If $n=1$, then $w_1$ is not the identity.
\end{numlist}

\end{defi}

The case $n=0$ is included as  the empty sequence. Also, if $n$ is larger than one then for each $i \in \{1,2,\ldots, 
n\}$, $w_{i} \notin C$, other wise, $(2)$ is violated.

The proof of the next theorem can be found in \cite{LS}, \cite{cohen} or \cite{MKS}.

\start{theo}{normal} \begin{numlist}
\item  Every element of $G \ast_C H$  can be written as a product $w_1w_2 \cdots w_n$ where $(w_1, w_2, \ldots, w_n)$ is a  reduced sequence.
\item If $n$ is a positive integer and $(w_1, w_2, \ldots, w_n)$ is a reduced sequence then the product $w_1w_2 \cdots w_n$ is not the identity.
\end{numlist}
\end{theo}

We could not find a proof of the next well known result in the literature, so we include it here.

\start{theo}{reduced amalgamation} If  $(w_1, w_2, \ldots, w_n)$  and
$(h_1, h_2, \ldots, h_n)$ are reduced sequences such the products
$w_1w_2\cdots w_n$ and $h_1 h_2, \cdots h_n$ are equal then exists a finite sequence of
element of C, namely $(c_1, c_2, \ldots, c_{n-1})$ such that
 \begin{romlist}
 \item $w_1=h_1 c_1$,
 \item $w_{n}=c_{n-1}^{-1} h_n$,
 \item For each  $i \in \{2,3,\ldots,n-1\}$, $w_i=c_{i-1}^{-1}h_i c_i$.
\end{romlist}
\end{theo}
\begin{proof}
We can assume $n>1$. Since the products are equal, $h_n^{-1}h_{n-1}^{-1}\cdots
h_1^{-1}w_1w_2\cdots w_n$ is the identity. By \theoc{normal}(2), the
sequence
$$(h_n^{-1},h_{n-1}^{-1},\ldots , h_1^{-1},w_1, w_2, \ldots, w_n)$$ is not reduced. Since the sequences ($h_n^{-1}, 
h_{n-1}^{-1}, \ldots, h_1^{-1})$ and  $(w_1, w_2,
\ldots, w_n)$   are reduced, $w_1$ and $h_1^{-1}$ belong both to
$G\setminus C$ or both to $H \setminus C$. Assume  the first possibility holds, that is $h_1$ and
$w_1$ are in  $G\setminus C$, (the second possibility is treated analogously.) Set $c_1= h_1^{-1}w_1$.
By our assumption, $c_1 \in G$. The sequence
\begin{equation}\nonumber
(h_n^{-1},h_{n-1}^{-1},\ldots, h_2^{-1}, c_1, w_2, \ldots, w_n)
\end{equation}
 has product equal to the identity. All the elements of this sequence are
alternatively in $G\setminus C$ and $H\setminus C$, with the possible exception of $c_1$. By
\theoc{normal}(2), this sequence is not reduced. Then $c_1 \notin G \setminus C$.  Since $c_1 \in
G$, $c_1 \in C$.

Now, consider the sequence
\begin{equation}\nonumber
(h_n^{-1},h_{n-1}^{-1},\ldots,h_3^{-1}, h_2^{-1} c_1 w_2,w_3, \ldots, w_n)
\end{equation}
By arguments analogous to those we made before, $h_2^{-1} c_1 w_2 \in C$. Denote
$c_2=h_2^{-1}c_1w_2$. Thus, $w_2=c_1^{-1}h_2c_2$

This shows that we can apply induction   to find  the sequence $(c_1, c_2, \ldots, c_{n-1})$ claimed
in the theorem.

\end{proof}

\start{defi}{cyclically reduced form} A finite sequence of elements $(w_1, w_2, \ldots,w_n)$ of $G
\ast_C H$ is  \emph{cyclically reduced } if every cyclic permutation of $(w_1, w_2,\ldots, w_n)$ is
reduced.
\end{defi}

\start{nota}{index} When dealing with free products with
amalgamation, every time we consider a sequence of the form $(a_1,
a_2, \ldots, a_n)$ we take subindexes mod $n$ in the following way:
For each $j \in \Z$, by $a_j$ we will denote $a_i$ where $i$ is the
only integer in $\{1,2,\ldots,n\}$ such that $n$ divides $i-j$.
\end{nota}

\start{rem}{even} If  $(w_1, w_2, \ldots, w_n)$ is a cyclically reduced
sequence and $n \ne 1$ then $n$ is even. Also, for every pair of integers $i$ and $j$, $w_i$ and $w_j$ are both in 
$G$ or both in $H$ if and only if $i$ and $j$ have the same parity. 
\end{rem}

The following result is a direct consequence of \theoc{normal}(1)  and  \cite[Theorem 4.6]{MKS}. 

\start{theo}{existence}  Let $s$ be a conjugacy class of  $G \ast_C H$. The there exists a cyclically reduced 
sequence  such that the product is a representative of  $s$. Moreover  every cyclically reduced sequence with 
product in $s$ has the same number of terms.
 \end{theo}

The following result gives necessary conditions for two cyclically reduced sequences to be conjugate.

 \start{theo}{conjugacy reduced} Let $n \ge 2$ and let  $(w_1, w_2,\ldots,w_n)$ and $(v_1, v_2, \ldots, v_n)$ be
cyclically reduced sequences such that the products $w_1w_2,\cdots w_n$ and $v_1 v_2 \cdots v_n$ are
conjugate. Then there exists an  integer $k \in \{0,1,\ldots,n-1\}$ and a sequence of elements of
$C$, $c_1, c_2, \ldots, c_n$ such that for each $i \in \{1,2,\ldots,n\}$,
$w_i=c_{i+k-1}^{-1}v_{i+k}c_{i+k}.$  In particular, for each $i \in \{1,2,\ldots,n\}$, $w_i$ and $v_{i+k}$ are both in $G
$ or both in $H$.
\end{theo}
\begin{proof}
By \cite[Theorem 2.8]{LS}, there exist $k \in \{0,1,2,\ldots, n-1\}$ and an element
$c$ in the amalgamating group $C$ such that $$w_1w_2\cdots w_n=c^{-1}v_{k+1}v_{k+2}\cdots
v_{k+n-1}v_{k+n}c.$$ The sequences $(w_1, w_2,\ldots,w_n)$ and $(c^{-1}v_{k+1},v_{k+2},\ldots,
v_{k+n-1},v_{k+n}c)$ are reduced  and have the same product. By \theoc{reduced amalgamation} there exists a 
sequence
of elements of $C$, $(c_1, c_2,\ldots, c_{n-1})$ such that
\begin{romlist}
 \item $w_1=c^{-1}v_{k+1} c_1$,
 \item $w_{n}=c_{n-1}^{-1} v_{k+n} c$,
 \item For each  $i \in \{2,3,\ldots,n-1\}$, $w_i=c_{i-1}^{-1}v_{k+i} c_i$.
\end{romlist}

Relabeling the sequence $(c, c_1, c_2, \ldots, c_{n-1})$ we obtain the desired result.
\end{proof}

If $C_1$ and $C_2$  are subgroups of a group $G$ a \emph{double coset of $G$ mod $C_1$ on the left and $C_2$ 
on the right } or briefly a \emph{double coset of $G$} is an equivalence class of the equivalence relation on $G$ 
defined  for each pair of elements $x$ and $y$  of $ G$ by $x \sim y$ if there exist $c_1 \in C_1$ and $c_2 \in C_2$ 
such that $x = c_1yc_2$. If $x \in G$, then the equivalence class containing $x$ is denoted by $C_1xC_2$.
 Using this notation, the next corollary follows directly from \theoc{conjugacy reduced}.

\start{cory}{double cosets} Let $n \ge 2$ and let  $(w_1, w_2,\ldots,w_n)$   and $(v_1, v_2, \ldots, v_n)$ be
cyclically reduced sequences such that the products $w_1w_2,\cdots w_n$ and $v_1 v_2 \cdots v_n$ are
conjugate. Then $$\{Cw_1w_2C, Cw_2w_3C,\ldots, Cw_nw_1C\}=\{Cv_1v_2C, Cv_2v_3C,\ldots, Cv_nv_1C\}.$$
 
\end{cory}

\begin{rem} Observe that a result stronger than \coryc{double cosets} holds, namely, one can associate a  unique 
cyclic sequence of double cosets (and not only a set) to each conjugacy class. \end{rem}

\start{defi}{malnormal subgroup} Let $C$ be a subgroup of a group  $G$ and let $g$ be an element of  $G$. 
Denote by $C^g$ the subgroup of $G$ defined by $g^{-1}Cg$. We say that \emph{
 $C$ is malnormal in $G$} if $C^g \cap C = \{1\}$ for every $g \in G\setminus C$.
\end{defi}

\start{lem}{claim1} Let $G \ast_C  H$ be an free product with amalgamation such that the amalgamating group $C$ 
is malnormal in $G$ and is malnormal in $H$. Let $a$ and $b$ be elements of  $C$. Let $w_1,w_2$ and $v_1,v_2
$  be two reduced sequences such that  the sets of double cosets
$$\{Cw_1aw_2C, Cv_1v_2C\} \mbox{ and }  \{Cw_1w_2C,Cv_1bv_2C\}$$ are equal. Then $a$ and $b$ are 
conjugate in $C$. Moreover, if $a\ne 1$ or $b \ne 1$ then $v_1$ and $w_1$ are both in $G$ or both in $H$. 
\end{lem}
\begin{proof}

We claim that  if $Cw_1a w_2 C=Cw_1 w_2 C$ then $a= \unitC$. Indeed, if  $Cw_1aw_2 C=Cw_1 w_2 C$ then  
there exist $c_1$ and $c_2$ in $C$ such that  $c_1 w_1 a w_2 c_2 =  w_1 w_2.$ By \theoc{reduced 
amalgamation} there exists $d \in C$ such that 
 \begin{align} 
 \label{ee2}
 w_1&=  c_1w_1 a d \\
  w_2& = d^{-1} w_2c_2 \label{ee1}
\end{align}

Since $C$ is malnormal in $H$, by Equation (\ref{ee1}), $d=1$. Thus Equation (\ref{ee2}) becomes $w_1=  c_1w_1 
a $. By malnormality, $a=1$ and the proof of the claim is complete.

Assume first that $Cw_1aw_2C= Cw_1w_2C$ and  $Cv_1v_2C =Cv_1bv_2C$ then by  the claim, $a=1$ and $b=1
$. Thus, $a$ and $b$ are conjugate in $C$ and the result holds in this case.  

Now, assume that $Cw_1aw_2C=Cv_1bv_2C$ and  $Cv_1v_2C = Cw_1w_2C$. By definition of double cosets, 
there exist  $c_1, c_2, c_3, c_4 \in C$ such that
$v_{1}bv_{2}=c_{1}w_{1}aw_{2}c_{2}$ and $v_{1}v_{2}=c_{3}w_{1}w_{2}c_{4}$.  By \theoc{reduced 
amalgamation},   there exist $d_1, d_2 \in C$ such that 
 \begin{align} 
 v_1 b&= c_1 w_1 a  d_1 \label{e3}\\
v_2 &=  d_1^{-1} w_2 c_2  \label{e4}\\
v_1  &= c_3 w_1   d_2 \label{e5}\\
v_2 &=  d_2^{-1} w_2 c_4  \label{e6}
\end{align}
By Equations (\ref{e4}) and (\ref{e6}) and malnormality, $d_1=d_2$. By By Equations (\ref{e3}) and (\ref{e5}) and 
malnormality, $d_1 =a d_1 b^{-1}$. Thus, $a$ and $b$ are conjugate in $C$ . Finally by Equation (\ref{e3}), since 
$C$ is a subgroup of $G$ and a subgroup of $H$, $w_1$ and $v_1$ are both in $G$ or both in $H$.

\end{proof}

Now we will show that certain pairs of conjugacy classes are distinct by proving that the associated sequences of 
double cosets of some of their respective representatives are distinct. Warren Dicks suggested the idea of this 
proof. Our  initial proof \cite{c1} applied repeatedly  the equations given by \theoc{conjugacy reduced} to derive a 
contradiction.

\start{theo}{amalgamating with cyclic}  Let  $G \ast_C H$ be a free
product with amalgamation. Let  $i$ and $j$ be distinct elements of  $\{1,2,\ldots, n\}$ and let $a$ and $b$ be 
elements of $C$. Assume all the following:
\begin{numlist}
\item  The subgroup $C$ is malnormal in $G$ and  is malnormal in $H$.
\item  If $i$ and $j$ have the same parity then $a$ and $b$ are not conjugate in $C$.
\item Either $a\ne 1$ or $b \ne 1$.
\end{numlist}
Then for every cyclically reduced sequence  $(w_1, w_2, \ldots w_n)$   the products
$$w_1w_2 \cdots w_i aw_{i+1}\cdots w_n \mbox{ and }  w_1w_2 \cdots w_j b w_{j+1}\cdots w_n  $$
are not conjugate. 
\end{theo}

\begin{proof}   Let $(w_1, w_2,\ldots, w_n)$ be a cyclically reduced sequence. Assume that there exist
$a$ and $b$ in $C$ and $i, j \in \{1,2,\ldots,n\}$  as in the hypothesis of the theorem  such that  the products
$w_1w_2 \cdots w_i aw_{i+1}\cdots w_n $ and $w_1w_2 \cdots w_j bw_{j+1}\cdots w_n $ are conjugate.

By \coryc{double cosets}, the sequences of double cosets mod $C$ on the right and the left associated with $(w_1, 
w_2,\ldots, aw_{i+1},\ldots ,w_n)$ and  with $(w_1, w_2,\ldots, bw_{j+1},\ldots ,w_n)$ are equal. In symbols,
\begin{equation}\nonumber
\{Cw_1w_2C, \ldots,Cw_iaw_{i+1}C, \ldots, Cw_nw_1C\}=\{Cw_1w_2C,\ldots, Cw_jb w_{j+1}C,\ldots, Cw_nw_1C\}
\end{equation}
Removing from both sets the elements that are denoted by equal expressions we obtain 
\begin{equation}\nonumber
\{Cw_iaw_{i+1}C, Cw_jw_{j+1}C\}=\{ Cw_jb w_{j+1}C, Cw_iw_{i+1}C\}.\end{equation}
By \lemc{claim1}, $a$ and $b$ are conjugate in $C$. Moreover, $w_i$ and $w_j$ are both in $G$ or both in $H$. 
By \remc{even}, $i$ and $j$ have the same parity,
 contradicting our hypothesis $(2)$.
\end{proof}

\start{rem}{hyp} It is not hard to construct an example that shows that hypothesis 
$(1)$ of \theoc{amalgamating with cyclic} is necessary. Indeed take, for instance, $G$ and $H$ two infinite cyclic 
groups generated by $x$ and $y$ respectively. Let $C$ be an infinite cyclic subgroup generated by $z$. Define $
\map{\varphi}{C}{G}$ and  $\map{\psi}{C}{H}$ by $\varphi(z)=x^{2}$ and $\psi(z)=y^{3}$. The sequence $(x,y)$ is 
cyclically reduced. Let $a=x^{2}$ and $b=x^{2}$. 
The products $x a y$ and $x y b$ are conjugate.

The following example shows that  hypothesis 
$(2)$ of \theoc{amalgamating with cyclic} is necessary. 
Let $G \ast_C H$ be a free product with amalgamation, let $c$ be an element of $C$ and let $(w_{1},w_{2})$ be a 
reduced sequence. Thus, $(w_{1},w_{2},w_{1},w_{2})$ is a cyclically reduced sequence. On the other hand, the 
products $w_{1} c  w_{2}w_{1}w_{2}$ and $w_{1}w_{2}w_{1} c w_{2}$ are conjugate. 

\end{rem}

The next result states certain elements of an amalgamated free product are not conjugate.

\start{theo}{amalgamating inverse}  Let $G \ast_C H$ be a free product with amalgamation. Let $a$ and $b$ be elements of
$C$ and let and  $i, j \in \{1,2,\ldots, n\}$.  Assume that for  every $g \in (G \cup H) \setminus C$,  are not in the same double coset relative to $C$, i.e., $g^{-1} \notin CgC$.
Then for every cyclically reduced sequence  $(w_1, w_2, \ldots w_n)$   the products
$$w_1w_2 \cdots w_i aw_{i+1}\cdots w_n \mbox{ and }  w_n^{-1}w_{n-1}^{-1} \cdots w_{j+1}^{-1} b w_{j}^{-1}\cdots 
w_1^{-1}  $$
are not conjugate.
\end{theo}
\begin{proof} Assume that the two products are conjugate. Observe that the sequences $$(w_1, w_2,  \ldots,  w_i 
a,  w_{i+1},\ldots, w_n) \mbox{ and }  (w_n^{-1}, w_{n-1}^{-1},  \ldots,  w_{j+1}^{-1} b,  w_{j}^{-1}, \ldots , w_1^{-1} )
$$ are cyclically reduced. By \theoc{conjugacy reduced}   there exists an
integer $k$  such that for every $h \in \{1,2,\ldots,n\}$ $w_h$ and $w_{1-h+k}^{-1}$ are both in $G$ or both in $H$. 
Then $w_h$ and $w_{1-h+k}$ are both in $G$ or both in $H$. By \remc{even}, $n$ is even and  $h$ and $1-h+k$ 
have the same parity. Thus, $k$ is odd.

Set $l=\frac{k+1}{2}$. By \theoc{conjugacy reduced},  $w_l$ and $w_{1-l+k}^{-1}=w_{l}^{-1}$  are in the same 
double coset of $G$ (if $l \in \{i, j-k\}$, either $w_l$ or $w_{1-l+k}$ may appear multiplied by $a$ or $b$ in the 
equations of \theoc{conjugacy reduced} but this does not change the double coset.) Thus $w_l$ and $w_{l}^{-1}$ 
are in the same double coset  mod $C$ of $G$, contradicting our hypothesis.

\end{proof}

\section{Oriented separating simple loops}\label{separating oriented}

The goal of this section is to prove \theoc{bracket form}, which gives us a way to compute the bracket of a 
separating simple closed curve $x$ and the product of the terms of a cyclically reduced sequence given by the 
amalgamated free product that $x$ determines.  We prove this theorem by finding in \lemc{separating 
representative} appropriate representatives for the separating simple closed curve $x$ and the  terms of the 
cyclically reduced sequence.

Through the rest of these pages, a \emph{surface} will mean a connected oriented surface. We denote
such a surface by $\Sigma$. The fundamental group of   $\Sigma$ will be denoted by
$\pi_1(\Sigma,p)$ where $p \in \Sigma$ is the basepoint. By a \emph{curve} we will mean a closed
oriented curve on $\Sigma$. We will use the same letter to denote   a curve and its image on $\Sigma$.

Let $\chi$ be a separating non-trivial simple curve on $\Sigma$, non-parallel to a boundary
component of $\Sigma$. Choose a point $p \in \chi$ to be the basepoint of each of  the fundamental
groups which will appear in this context. Denote by $\Sigma_1$ the union of $\chi$ and one of the
connected components of $\Sigma \setminus \chi$ and by $\Sigma _2$ the union of $\chi$ with the
other connected component.

\start{rem}{canonical} As a consequence of the van Kampen's theorem (see \cite{h})
$\pi_1(\Sigma,p)$ is canonically isomorphic to the free product of $\pi_1(\Sigma_1,p)$ and
$\pi_1(\Sigma_2,p)$ amalgamating the subgroup $\pi_1(\chi,p)$, where the monomorphisms $\pi_1(\chi,
p)\longrightarrow \pi_1(\Sigma_1,p)$ and $\pi_1(\chi, p)\longrightarrow \pi_1(\Sigma_2,p)$ are the
induced by the respective inclusions.
\end{rem}

\start{lem}{separating representative} Let $\chi$ be a separating simple closed curve on $\Sigma$.
Let $(w_1, w_2, \ldots ,w_n)$ be a cyclically reduced sequence for the amalgamated product of
\remc{canonical}. Then there exists a sequence of curves $(\gamma_1, \gamma_2, \ldots, \gamma_n)$ such
that for each $i \in \{1,2,\ldots,n\}$ the all the following holds.
\begin{numlist}
\item  The curve $\gamma_i$ is a representative of $w_i$.
\item  The curve $\gamma_i$ is alternately in $\Sigma_1$ and $\Sigma_2$.
\item  The point $p$ is the basepoint of $\gamma_i$.
\item The point  $p$ is not a self-intersection point of $\gamma_i$. In other words, $\gamma_i$ passes through
$p$ exactly once.
\end{numlist}
Moreover, the product $\gamma_{1}\gamma_{2}\cdots \gamma_{n}$ is a representative of the product $w_{1}w_
{2}\cdots w_{n}$ and the curve $\gamma_{1}\gamma_{2}\cdots \gamma_{n}$  and $w_{1}w_{2}\cdots w_{n}$  
intersect $\chi$  transversally with multiplicity $n$. 
\end{lem}
\begin{proof} For each $i \in \{1,2,\ldots,n\}$, take  a representative $\gamma_i$ in of $w_i$. Notice that  $
\gamma_i \subset\Sigma_1$ or $\gamma_i\subset \Sigma_2$. We   homotope $\gamma_i$ if
necessary, in such a way that $\gamma_i$ passes through $p$ only once, at the basepoint $p$

\begin{figure}[htbp]\label{lemma separating}
\begin{center}
\begin{pspicture}(12,5)
 \psccurve(1,2)(1,0.2)(6,1)(9,0.3)(11,3)(6,3)(1,4.4)
\pscurve(2,1.5)(2.7,1.2)(3.5,1.5)
\pscurve(2.2,1.4)(2.7,1.5)(3.2,1.34)

\rput(-0.5,1.58){\pscurve(2,1.5)(2.7,1.2)(3.5,1.5)
\pscurve(2.2,1.4)(2.7,1.5)(3.2,1.34) }
\rput(6,0.51){\pscurve(2,1.5)(2.7,1.2)(3.5,1.5)
\pscurve(2.2,1.4)(2.7,1.5)(3.2,1.34) }

\pscurve[linecolor=gray,linewidth=1pt](6,1)(6.15,2)(6,3)
\pscurve[linecolor=gray,linewidth=1pt,linestyle=dashed](6,1)(5.9,2)(6,3)

\psccurve[showpoints=false,linecolor=lightgray](1,3)(3,2.5)(6.15,2.2)(9.2,1.5)(9.6,3)(6.15,2.2)(2,4)(2,1)(6.15,2.2)
(10,1.9)(9,2.5)(6.15,2.2)

\rput(6.2,3.3){$\chi$}

\rput(10,3){$\gamma_1$}
\rput(10.3,2){$\gamma_3$}
\rput(1,3.4){$\gamma_2$}
\rput(1.8,4){$\gamma_4$}

\end{pspicture}

\end{center}
\caption{The representative of Lemma \ref{separating representative} }
\end{figure}
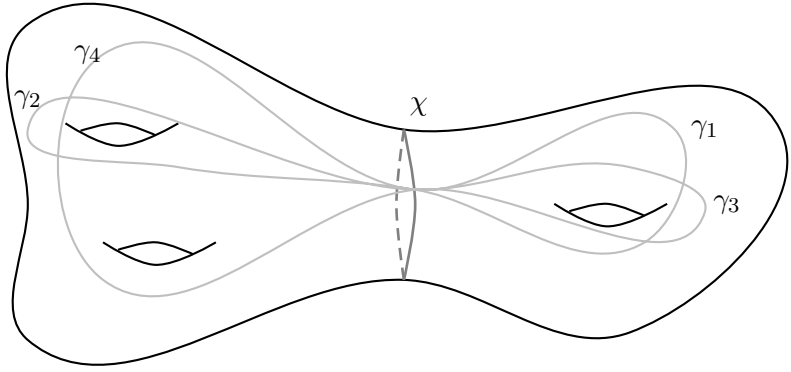

Since each $\gamma_{i}$ intersects $\chi$ only at $p$, the product  $\gamma_{1}\gamma_{2}\cdots \gamma_{n}$  
intersects $\chi$ exactly $n$ times. Each of this intersections happens when the curve $\gamma_{1}\gamma_{2}
\cdots \gamma_{n}$ passes from $\Sigma_{1}$ to $\Sigma_{2}$ or from $\Sigma_{2}$ to $\Sigma_{1}$. This 
implies that these $n$ intersection points of $\chi$ with $\gamma_{1}\gamma_{2}\cdots \gamma_{n}$ are 
transversal.  

\end{proof}

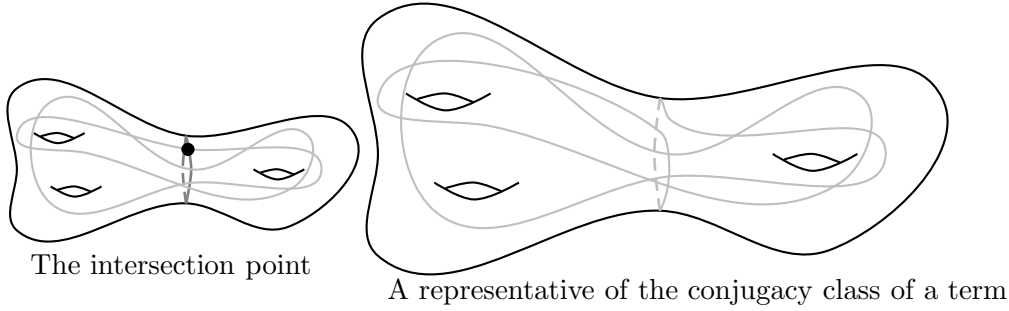
\begin{figure}[ht]\label{term separating}
\begin{center}
\begin{pspicture}(12,3.5)
\rput(-0.5,0.4){
\psset{xunit=0.45,yunit=0.45}
 \psccurve(1,2)(1,0.2)(6,1)(9,0.3)(11,3)(6,3)(1,4.4)
\pscurve(2,1.5)(2.7,1.2)(3.5,1.5)
\pscurve(2.2,1.4)(2.7,1.5)(3.2,1.34)

\rput(-0.5,1.58){\pscurve(2,1.5)(2.7,1.2)(3.5,1.5)
\pscurve(2.2,1.4)(2.7,1.5)(3.2,1.34) }
\rput(6,0.51){\pscurve(2,1.5)(2.7,1.2)(3.5,1.5)
\pscurve(2.2,1.4)(2.7,1.5)(3.2,1.34) }

\pscurve[linecolor=gray,linewidth=1pt](6,1)(6.15,2)(6,3)
\pscurve[linecolor=gray,linewidth=1pt,linestyle=dashed](6,1)(5.9,2)(6,3)

\psccurve[showpoints=false,linecolor=lightgray](1,3)(3,2.5)(9.2,1.5)(9.6,3)(6.2,2)(2,4)(2,1)(6.1,1.6)(10,2)(6.1,2.6)

}

\psdot*[dotsize=5pt](2.225,1.566)

\rput(9,-0.320){A representative of the  conjugacy class of a term}

\rput(2,0){The intersection point}

\rput(4,0){
\psset{xunit=0.75,yunit=0.75}

 \psccurve(1,2)(1,0.2)(6,1)(9,0.3)(11,3)(6,3)(1,4.4)
\pscurve(2,1.5)(2.7,1.2)(3.5,1.5)
\pscurve(2.2,1.4)(2.7,1.5)(3.2,1.34)

\rput(-0.5,1.58){\pscurve(2,1.5)(2.7,1.2)(3.5,1.5)
\pscurve(2.2,1.4)(2.7,1.5)(3.2,1.34) }
\rput(6,0.51){\pscurve(2,1.5)(2.7,1.2)(3.5,1.5)
\pscurve(2.2,1.4)(2.7,1.5)(3.2,1.34) }
\pscurve[linecolor=lightgray,linewidth=1pt,linestyle=dashed](6,1)(5.9,2)(6,3)

\psecurve[showpoints=false,linecolor=lightgray](5,0)(6,1)(6.15,2)(6,2.4)(1,3)(3,2.5)(9.2,1.5)(9.6,3)(6.2,2)(2,4)(2,1)
(6.1,1.6)(10,2)(6.2,2.6)(6,3)(5,4)
}

\end{pspicture}
\end{center}
\caption{An intersection point (left) and the corresponding term of the bracket (right)}
\end{figure}

Let $\Sigma$ be an oriented surface.
The Goldman bracket \cite{g} is a Lie bracket defined on the  vector space  generated by all free homotopy classes 
of  oriented curves on the surface $\Sigma$.  
We recall the definition: For each pair of homotopy classes $a$ and $b$, consider representatives $a$ and $b$ 
respectively, that only intersect  in transversal double points. The bracket of $[a, b]$ is defined as the signed sum 
over all intersection points $P$  of $a$ and $b$ of free homotopy class of the curve that goes around $a$ starting 
and
ending at $P$ and then goes around $b$ starting and ending at $P$. The sign of the term at an intersection point 
$P$ is the intersection number of $a$ and $b$ at $P$. (See \figc{term separating}.)

The above definition can be extended to consider pairs of representatives where branches intersect transversally but triple (and higher) points are allowed. Indeed, take a pair of such representatives $a$ and $b$. The Goldman Lie 
bracket is the sum over the intersection of pairs of small arcs, of the conjugacy classes of the curve obtained by 
starting in an intersection point and going along the $a$ starting in the direction of the first arc, and then going 
around $b$ starting in the direction of the second arc. The sign is determined by the pair of tangents of the  
ordered oriented arcs  at the intersection point.

\start{nota}{declaration} The Goldman bracket $[a,b]$ is
computed for pairs of conjugacy classes $a$ and $b$ of curves on
$\Sigma$. In order to make the notation lighter, we will abuse
notation by writing $[u,v]$, where $u$, $v$ are elements of the
fundamental group of $\Sigma$. By $[u,v]$, then, we will mean the
bracket of the conjugacy class of $u$ and the conjugacy class of
$v$.
\end{nota}

\start{theo}{bracket form} Let $x$  be a conjugacy class of $\pi_1(\Sigma,p)$ which can be
represented by a separating simple closed curve $\chi$. Let $(w_1, w_2, \ldots ,w_n)$ be a cyclically
reduced sequence of the for the amalgamated product determined by $\chi$ in \remc{canonical}.  Then $[w_1,x]=0
$. 
Moreover, if $n >1$ then there exists $s \in \{1,-1\}$ such that the bracket is given by
$$s  [w_1w_2\cdots w_n,x]= \sum_{i=1}^n (-1)^{i}w_1w_2\cdots w_i
xw_{i+1}\cdots w_n.$$
\end{theo}

\begin{proof} For each curve $\gamma$ there exists a representative of the null-homotopic class disjoint from $
\gamma$ so the result holds when $n=0$.
 We prove now  that $[w_1,x]=0$. By definition of a cyclically reduced sequence, $w_1 \in
\pi_1(\Sigma_1,p)$ or $w_1 \in \pi_1(\Sigma_2,p)$. Suppose that $w_1 \in \pi_1(\Sigma_1,p)$ (the
other possibility is analogous.) Choose
 a curve $\gamma_1 \subset \Sigma_1$ which is a representative of $w_1$. Since $\gamma_1 \subset \Sigma_1$ , 
we can
homotope $\gamma_1$ to a curve  which has no intersection with
$\chi$, (clearly, this is a free homotopy  which does not fix the
basepoint $p$.) This shows that $w_1$ and $x$ have disjoint
representatives, and then $[w_1,x]=0$.

Now, assume that $n>1$. Let $(\gamma_1, \gamma_2, \ldots, \gamma_n)$ be the sequence given by \lemc
{separating
representative} for $(w_1, w_2, \ldots ,w_n)$ .

The loop product $\gamma_1\gamma_2\cdots\gamma_n$ is a representative of the group product
$w_1w_2\cdots w_n$.

Every intersection  of $\gamma_1\gamma_2\cdots\gamma_n$ with $\chi$ occurs when
$\gamma_1\gamma_2\cdots\gamma_n$ leaves  one connected component of $\Sigma \setminus \chi$ to enter
the other connected component. That is, between each  $\gamma_i$ and $\gamma_{i+1}$. (Recall we are using
\notc{index} so the intersection between $\gamma_n$ and $\gamma_1$ is considered.)

For each $i \in \{1,2,\ldots,n\}$ denote by $p_i$ the intersection point of $\chi$ and
$\gamma_1\gamma_2\cdots \gamma_n$ between $\gamma_i$ and $\gamma_{i+1}$. The loop product 
$\gamma_1\gamma_2\cdot\gamma_i \chi \gamma_{i+1} \cdots\gamma_n$ is a representative of the conjugacy 
class of the term of the Goldman Lie bracket corresponding to $p_i$. Thus, for each $i \in \{1,2,\ldots,n\}$ , the 
conjugacy class of the term of the bracket
corresponding to the intersection point  $p_i$ has $w_1\cdots w_i x w_{i+1}\cdots w_n$ as representative.

Let $i, j \in \{1,2,\ldots,n\}$ with different parity. Assume that $w_i \in \pi_1(\Sigma_1,p)$, (the
case $w_i \in \pi_2(\Sigma_1,p)$ is similar.) The tangent vector of $\gamma_1\gamma_2\cdots
\gamma_n$ at $p_i$ points towards $\Sigma_2$ and the tangent vector of $\gamma_1\gamma_2\cdots
\gamma_n$ at $p_j$ points towards $\Sigma_1$. This shows that the signs of the bracket terms
corresponding to $p_i$ and $p_j$ are opposite, completing  the proof.
\end{proof}

\start{rem}{triple} In \theoc{bracket  form}, all the intersections of the chosen representatives
of $w_1w_2\cdots w_n$ and $x$ occur at the basepoint $p$. The representative $\omega$ of
$w_1w_2\cdots w_n$ intersects $x$ in a point which  is a multiple self-intersection point of $\omega$.
This does not present any difficulty in computing the bracket, because the intersection points of
both curves are still transversal double points.
\end{rem}

\section{HNN extensions}\label{hnn extensions}

This section is the HNN counterpart of  Section \ref{Amalgamated products} and the statements, arguments and 
posterior use of the statements are similar.  The main goal consists in proving that the products of certain 
cyclically reduced sequences cannot be conjugate (Theorems \ref{non-separating hnn} and \ref{hnn inverse}.) This 
result will be used to show that the pairs of terms of the bracket of certain conjugacy classes with opposite sign do 
not cancel.

Let $G$ be a group, let $A$ and $B$ be two subgroups of $G$ and let $\map{\varphi}{A}{B}$ be an
isomorphism. Then the \emph{HNN extension of $G$ relative to $A,B$ and $\varphi$ with stable letter
$t$ (or, more briefly, the HNN extension of $G$ relative to $\varphi$)} will be denoted by $G^{\ast
\varphi}$ and is the group obtained by taking the quotient of the free product of $G$ and the free group generated 
by $t$ by the normal subgroup generated by $t^{-1}at\varphi(a)^{-1}$ for all $a \in A$.
 (see \cite{LS} for detailed definitions.)


\start{defi}{reduced hnn}  Consider an HNN extension  $G^{\ast \varphi}$. Let $n$ be a non-negative integer and 
for each $i \in \{1,2,\ldots,n\}$, let $\varepsilon_i\in \{1,-1\}$ and $g_i$ be an element of $G$. A finite sequence $
(g_0,t^{\varepsilon_1},g_1,t^{\varepsilon_2}\ldots,
g_{n-1},t^{\varepsilon_n},g_n)$  is said to be \emph{reduced} if there is no consecutive subsequence
of the form $(t^{-1},g_i, t)$ with $g_i \in A$ or $(t, g_j, t^{-1})$ with $g_j \in B$.
\end{defi}

The following result is the analogue of \theoc{normal} for HNN extensions (see \cite{LS} or \cite{cohen}.)

\start{theo}{Britton}\begin{numlist}
\item (Britton's lemma) If the sequence
$(g_0,t^{\varepsilon_1},g_1,t^{\varepsilon_2}\ldots, g_{n-1},t^{\varepsilon_n},g_n)$ is reduced and $n
\ge 1$ then the product  $g_0t^{\varepsilon_1}g_1t^{\varepsilon_2}\cdots
g_{n-1}t^{\varepsilon_n}g_n$ is not the identity in the HNN extension $G^{\ast \varphi}$.

\item Every element  $g$  of $G^{\ast \varphi}$ can be written as a product
$g_0t^{\varepsilon_1}g_1t^{\varepsilon_2}\cdots g_{n-1}t^{\varepsilon_n}g_n$ where the sequence
$(g_0,t^{\varepsilon_1},g_1,t^{\varepsilon_2}\ldots, g_{n-1},t^{\varepsilon_n},g_n)$ is  reduced.

\end{numlist}
\end{theo}

As in the case of  \theoc{reduced amalgamation}, we include the proof of the next known  result because we were 
unable to find it in the literature.

\start{theo}{main hnn} Suppose that the equality
$$g_0t^{\varepsilon_1}g_1t^{\varepsilon_2}\cdots
g_{n-1}t^{\varepsilon_n}g_n=h_0t^{\eta_1}h_1t^{\eta_2}\cdots h_{n-1}t^{\eta_n}h_n$$ holds where
$(g_0,t^{\varepsilon_1},g_1,t^{\varepsilon_2}\ldots, g_{n-1},t^{\varepsilon_n},g_n)$ and
$(h_0,t^{\eta_1},h_1,t^{\eta_2}\ldots, h_{n-1},t^{\eta_n},h_n)$ are reduced sequences.

Then for each $i \in \{1,2,\ldots,n\}$ we have that $\varepsilon_i=\eta_i$. Moreover,  there exists a
sequence of elements $(c_1, c_2, \ldots, c_n)$ in $A \cup B$ such that
\begin{numlist}
\item $g_0 =h_0c_1  $
\item $g_n =\varphi^{\varepsilon_n}(c_{n}^{-1})h_n $
\item For each $i \in \{1,2,\ldots,n-1\}$,
$g_i=\varphi^{\varepsilon_i}(c_{i}^{-1})h_i c_{i+1}$
\item For each $i \in \{1,2,\ldots,n\}$,   $c_i \in A$ if
$\varepsilon_i=1$ and $c_i \in B$ if $\varepsilon_i=-1$.
\end{numlist}
\end{theo}
\begin{proof}

If the two products are equal then
\begin{equation}\label{q}
h_n^{-1}t^{-\eta_n} \cdots t^{-\eta_2}h_1^{-1}t^{-\eta_1}h_0^{-1} g_0t^{\varepsilon_1}g_1
t^{\varepsilon_2}\cdots g_{n-1}t^{\varepsilon_n}g_n =1
\end{equation}

By  \theoc{Britton}(1) the sequence that yields the  product on the left hand side of Equation
(\ref{q}) is not reduced. This implies that $\varepsilon_1 = \eta_1$. Moreover, if
$\varepsilon_1=1$ then
  $h_0^{-1} g_0 \in A$ and if $\varepsilon_1=-1$ then $h_0^{-1} g_0 \in B$.

Denote  the product $h_0^{-1}g_0$ by $c_1$. Thus, $g_0=h_0c_1$. By definition of HNN extension, we
can replace $t^{-\varepsilon_1}c_1 t^{\varepsilon_1}$ by $\varphi^{\varepsilon_1}(c_1)$ in Equation
(\ref{q}) to obtain,
\begin{equation}\label{qq}
h_n^{-1}t^{-\eta_n} \cdots t^{-\eta_2}h_1^{-1} \varphi^{\varepsilon_1}(c_1) g_1
t^{\varepsilon_2}\cdots g_{n-1}t^{\varepsilon_n}g_n =1
\end{equation}

By \theoc{Britton}, the sequence yielding the product of the left hand side of Equation (\ref{qq})
is not reduced. Hence, $\varepsilon_2=\eta_2$ and  if $\varepsilon_2=1$ then $h_1^{-1}
\varphi^{\varepsilon_1}(c_1) g_1 \in A$ and if $\varepsilon_2=-1$ then $h_1^{-1}
\varphi^{\varepsilon_1}(c_1) g_1 \in B$.

Denote by $c_2$ the product $h_1^{-1} \varphi^{\varepsilon_1}(c_1) g_1$. Thus, $g_1=
\varphi^{\varepsilon_1}(c_1^{-1})   h_1 c_2$.

By applying these  arguments, we can complete the proof by induction.
\end{proof}

\start{defi}{cyclic hnn}  Let $n$ be a non-negative integer. A  sequence  of elements of $G^{\ast \varphi}$, $(g_0,t^
{\varepsilon_1},g_1,t^{\varepsilon_2}, \ldots, g_{n-1},t^{\varepsilon_n})$ is said to be \emph{cyclically reduced} if  
all its cyclic permutations of
 are reduced.
\end{defi}

We could not find a direct proof in the literature of the first statement  of our next so we include  it here. (The 
second statement also follows from \theoc{cyclic hnn form} but it is a direct consequence of our proof.)
\start{theo}{existence hnn} Let $s$ be a conjugacy class of  $G^{\ast
\varphi}$. Then there exists a cyclically reduced sequence  such that
the product of its terms  is a representative of  $s$. Moreover, every cyclically
reduced sequence with product in $s$ has the same number of terms.
\end{theo}
\begin{proof}
If $s$ has a representative in $G$,  the result follows directly. So we can assume that $s$ has no
representatives in $G$.

By \theoc{Britton}(2), the set of reduced sequences with product in
$s$ is not empty.  Thus, it is possible to choose among all such
sequences, one that makes the number of terms the smallest possible.
Let $(g_0,t^{\varepsilon_1},g_1,t^{\varepsilon_2}, \ldots,
g_{m-1},t^{\varepsilon_m})$ be such a sequence. We claim that
$(g_0,t^{\varepsilon_1},g_1,t^{\varepsilon_2}, \ldots,
g_{m-1},t^{\varepsilon_m})$ is cyclically reduced.

Indeed, if  $m \in \{0,1\}$, the sequence has the form $(g_{0})$ or
$(g_0,t^{\varepsilon_1})$ and so it  is cyclically reduced. Assume now that
$m>1$.

If $(g_0,t^{\varepsilon_1},g_1,t^{\varepsilon_2}, \ldots, g_{m-1},t^{\varepsilon_m})$ is not cyclically
reduced then one of the following statements holds: \begin{numlist}
\item  $\varepsilon_1=1$, $\varepsilon_m=-1$ and $g_0 \in A$.
\item $\varepsilon_1=-1$, $\varepsilon_m=1$ and $g_0 \in B$.
\end{numlist}

We prove the result in case $(1)$. (Case $(2)$ can be treated with similar ideas.) In this case,
$t^{-\varepsilon_m}g_0t^{\varepsilon_1}=\varphi(g_0) \in B$

The sequence $(g_{m-1}\varphi(g_0)g_1,t^{\varepsilon_2}, \ldots,g_{m-2}, t^{\varepsilon_{m-1}}) $ is
reduced,   has product in $s$ and has strictly fewer terms than the sequence
$(g_0,t^{\varepsilon_1},g_1,t^{\varepsilon_2}, \ldots, g_{m-1},t^{\varepsilon_m})$. This contradicts our
assumption that our original sequence has the smallest number of terms.
Thus, our proof is complete.
\end{proof}

\start{nota}{permutations hnn} By definition,
the cyclically reduced  sequences of given HNN extension have the form
$(g_0,t^{\varepsilon_1},g_1,t^{\varepsilon_2}, \ldots,$ $
g_{m-1},t^{\varepsilon_m})$. From now on, we will make use of the following convention: For every integer $h$, $g_h$ will
denote $g_i$ where $i$ is the unique integer in $\{0,1,2, \ldots,
n-1\}$ such that $n$ divides $i-h$. Analogously,
 $\varepsilon_h$ will denote $\varepsilon_i$ where $i$ is  the unique integer  in $\{1,2,\ldots, n\}$ such that $n$ divides $i-h$.  \end{nota}

The next result is due to Collins and  gives necessary  conditions for two cyclically reduced sequence have 
conjugate product(see \cite{LS}.) 

\start{theo}{cyclic hnn form} (Collins' Lemma) Let $n \ge 1$ and let
$(g_0,t^{\varepsilon_1},g_1,t^{\varepsilon_2}, \ldots,
g_{n-1},t^{\varepsilon_n})$ and $(h_0,t^{\eta_1},h_1,t^{\eta_2},$
$\ldots, h_{m-1},t^{\eta_m})$ be two cyclically reduced sequences such
that their products are conjugate. Then $n=m$ and there exist $c
\in A \cup B$ and  $k \in \{0,1,2,\ldots, n-1\}$ such that the
following holds:
\begin{numlist}
\item $\eta_{k}=\varepsilon_n$,
\item $c \in A$ if $\varepsilon_n=-1$ and $c \in B$ if $\varepsilon_n=1$,
\item $g_0t^{\varepsilon_1} g_1 t^{\varepsilon_2}\cdots  g_{n-1}t^{\varepsilon_n}=c^{-1} h_{k}
t^{\eta_{k+1}} h_{k+1} t^{\eta_{k+2}}\cdots  h_{k+n-1}
t^{\eta_{k}}c.$
\end{numlist}

\end{theo}

By Theorems \ref{main hnn} and \ref{cyclic hnn form} and arguments exactly like those of
\theoc{conjugacy reduced}, we obtain the following result.

\start{theo}{hnn cory} Let $(g_0,t^{\varepsilon_1},g_1,t^{\varepsilon_2}, \ldots,
g_{n-1},t^{\varepsilon_n})$ and $(h_0,t^{\eta_1},h_1,t^{\eta_2}, \ldots, h_{n-1},t^{\eta_n})$ be
cyclically reduced sequences such that the products
$$g_0t^{\varepsilon_1}g_1t^{\varepsilon_2}\cdots
g_{n-1}t^{\varepsilon_n} \mbox{ and } h_0t^{\eta_1}h_1t^{\eta_2}\cdots h_{n-1}t^{\eta_n}$$ are
conjugate. Moreover, assume that $n \ge 1$. Then there exists an integer $k$ such that  for each $i \in \{1,2,
\ldots,n\}$,
$\varepsilon_i=\eta_{i+k}$. Moreover,  there exists a sequence of elements $(c_1, c_2, \ldots, c_n)$ in
$A \cup B$ such that for each $i \in \{1,2,\ldots,n\}$, $c_i \in A$ if $\varepsilon_i=1$ and $c_i
\in B$ if $\varepsilon_i=-1$ and
$$
g_i=\varphi^{\eta{i+k}}(c_{i+k}^{-1})h_{i+k} c_{i+k+1}$$
\end{theo}

\start{nota}{C} Let  $G^{\ast\varphi}$ be an HNN extension, where $\map{\varphi}{A}{B}$. We denote the subgroup 
$A$ by $C_1$ and the subgroup $B$ by  $C_{-1}$.
\end{nota}

\start{cory}{hnn cory cory}  Let $(g_0,t^{\varepsilon_1},g_1,t^{\varepsilon_2}\ldots,
g_{n-1},t^{\varepsilon_n})$ and $(h_0,t^{\eta_1},h_1,t^{\eta_2}, \ldots, h_{n-1},t^{\eta_n})$ be
cyclically reduced sequences such that the products
$$g_0t^{\varepsilon_1}g_1t^{\varepsilon_2}\cdots
g_{n-1}t^{\varepsilon_n} \mbox{ and } h_0t^{\eta_1}h_1t^{\eta_2}\cdots h_{n-1}t^{\eta_n}$$ are
conjugate.  If $n \ge 1$ then there exists an integer $k$ such that 
 for each $i \in \{1,2,\ldots,n\}$, $\varepsilon_i=\eta_{i+k}$ and  $
g_i$  belongs to  the double coset $C_{-\varepsilon_{i}} h_{i+k} C_{\varepsilon_{i+1}}.$
\end{cory}

\start{rem}{epsilon} If $\varepsilon \in \{1,-1\}$ and $a \in C_{\varepsilon}$ then $a  t^{\varepsilon} = t^{\varepsilon} 
\varphi^{\varepsilon}(a)$ in $G^{\ast\varphi}$.
\end{rem}

\start{defi}{dc1}  Let $n \ge 2$ and let  $(g_0,t^{\varepsilon_1},g_1,t^{\varepsilon_2}, \ldots,
g_{n-1},t^{\varepsilon_n})$ be a cyclically reduced sequence. 
The
\emph{sequence of double cosets associated with   $(g_0,t^{\varepsilon_1},g_1,t^{\varepsilon_2}, \ldots,
g_{n-1},t^{\varepsilon_n})$} is the sequence of double cosets $$(C_{\varepsilon_i}t^{{\varepsilon_i}} g_i t^{\varepsilon_{i+1}}g_{i+1}
t^{\varepsilon_{i+2}}C_{-\varepsilon_{i+2}})_{0 \le i \le n-1}.$$
\end{defi}

\start{lem}{claim two} Let $(g_0,t^{\varepsilon_1},g_1,t^{\varepsilon_2}\ldots, g_{n-1},t^{\varepsilon_n})$ and $
(h_0, t^{\eta_1}, h_1,t^{\eta_2}, , \ldots h_{n-1}, t^{\eta_n})$ be two cyclically reduced sequences whose products 
are conjugate and such that $n \ge 2$. Then the sequence of double cosets associated with  $(g_0,t^
{\varepsilon_1},g_1,t^{\varepsilon_2}\ldots, g_{n-1},t^{\varepsilon_n})$ is a cyclic permutation of the cyclic of double 
cosets associated with $(h_0, t^{\eta_1}, h_1t^{\eta_2}, \ldots h_{n-1}, t^{\eta_n})$. 
\end{lem}
\begin{proof} 
Let $k \in \{1,2,\ldots, n\}$ and  $(c_1, c_2, \ldots, c_n)$ be a finite  sequence of elements in $A \cup B$ given  by 
\theoc{hnn cory}  for the sequences $(g_0,t^{\varepsilon_1},g_1,t^{\varepsilon_2}\ldots, g_{n-1},t^{\varepsilon_n})$ 
and $(h_0, t^{\eta_1}, h_1,t^{\eta_2}, \ldots h_{n-1}, t^{\eta_n})$. By cyclically rotating  $(h_0, t^{\eta_1}, h_1,t^{\eta_2}, \ldots h_{n-1}, t^{\eta_n})$ if necessary, we can assume that $k=0$. (We could carry out the proof with $k>0$, but the assumption $k=0$ makes the equations more neat).

Let $i \in \{1,2,\ldots,n\}$. We will complete this proof by showing  that the $i$-th double coset of  the sequence $
(g_0, t^{\varepsilon_1} ,g_1, t^{\varepsilon_2}\ldots, g_{n-1},t^{\varepsilon_n})$ equals the $i$-th double coset 
of the sequence of  $(h_0, t^{\eta_1}, h_1,t^{\eta_2}, \ldots h_{n-1}, t^{\eta_n})$.

By \theoc{hnn cory}, $\varepsilon_j=\eta_{j}$ for each $j \in \{1,2,\ldots,n\}$ and 
$$
g_i t^{\varepsilon_{i+1}}g_{i+1}=\varphi^{\eta_{i}}(c_{i}^{-1})h_{i} c_{i+1} t^{\eta_{i +1}}\varphi^{\eta_{i +1}}(c_{i +1}^{-1})h_{i +1} c_{i +2}$$
By \remc{epsilon} we have that $ c_{i +1} t^{\eta_{i +1}}\varphi^{\eta_{i +1}}(c_{i +1}^{-1})=t^{\eta_{i +1}}.$ Thus, 
$$
g_i t^{\varepsilon_{i+1}}g_{i+1} =\varphi^{\eta_{i }}(c_{i }^{-1})h_{i } t^{\eta_{i +1}}h_{i 
+1}c_{i +2}
$$
By \theoc{hnn cory},  $c_{i }\in C_{\varepsilon_{i }}$ and $c_{i +2}\in C_{\varepsilon_{i +2}}$. Therefore, $\varphi^
{\eta_{i }}(c_{i }^{-1}) \in C_{-\varepsilon_{i }}$. Thus,
$$
t^{\varepsilon_{i}}g_i t^{\varepsilon_{i+1}}g_{i+1}t^{\varepsilon_{i+2}} \in t^{\varepsilon_{i}}C_{-\varepsilon_{i}}h_{i } t^{\eta_{i +1}}h_{i +1}C_{\varepsilon_{i+2}}t^{\varepsilon_{i}}$$
Consequently, by \remc{epsilon}, $t^{\varepsilon_{i}}g_i t^{\varepsilon_{i+1}}g_{i+1}t^{\varepsilon_{i+2}}$ 
  is in the $i$-th double coset associated with $(h_0, t^{\eta_1}, h_1,t^{\eta_2}, \ldots h_{n-1}, t^{\eta_n})$ and our proof is complete.
\end{proof}

\start{defi}{separated} Let  $G^{\ast\varphi}$ be an HNN extension, where $\map{\varphi}{A}{B}$. We say that $G^{\ast\varphi}$  is \emph
{separated} if $A \cap B^g=\{1\}$ for all $g \in G$. \end{defi}

\start{lem}{claim three}  Let $(t^{ \varepsilon_{1}},g_1,t^{ \varepsilon_{2}},g_2,t^{\varepsilon_3})$ and $(t^{\eta_{1}},h_1, t^{\eta_{2}}, h_2,t^{\eta_3})$ be 
two reduced sequences. Let $a \in C_{\varepsilon_{2}}$ and let $b \in C_{{\eta_{2}}}$.  Suppose that $G^{\ast \varphi}$  is 
separated and that $A$ and $B$ are malnormal in $G^{\ast\varphi}$.  Then the following statements hold.
\begin{numlist}
\item If the  double cosets  $C_{\varepsilon_{1}}t^{ \varepsilon_{1}} g_1 at^{ \varepsilon_{2}} g_2 t^{\varepsilon_3}C_{-\varepsilon_{3}}$ and 
$C_{\varepsilon_{1}}t^{ \varepsilon_{1}} g_1 t^{ \varepsilon_{2}} g_2 t^{\varepsilon_3}C_{-\varepsilon_{3}}$ are equal  then $a=1$.
\item If the subsets of double cosets $\{C_{\varepsilon_{1}}t^{ \varepsilon_{1}} g_1a t^{ \varepsilon_{2}} g_2 t^{\varepsilon_3}C_{-\varepsilon_{3}}, 
C_{\eta_{1}}t^{\eta_{1}} h_1  t^{\eta_{2}}  h_2 t^{\eta_3}C_{\varepsilon_{3}}\}$ and $\{C_{\varepsilon_{1}}t^{ \varepsilon_{1}} g_1 t^{ \varepsilon_{2}} g_2 t^{\varepsilon_3}C_{-\varepsilon_{3}}, 
C_{\eta_{1}}t^{\eta_{1}} h_1b  t^{\eta_{2}}  h_2 t^{\eta_3}C_{\varepsilon_{3}}\}$
are equal then $a$ and $b$ are conjugate by an element of $A \cup B$.  Moreover, if $a \ne 1$ or $b 
\ne 1$, then $ \varepsilon_{2}=\eta_{2}$.
\end{numlist}

\end{lem}
\begin{proof} We first prove $(1)$. 
By \remc{epsilon},
$$t^{ \varepsilon_{1}}C_{-\varepsilon_{1}} g_1 t^{ \varepsilon_{2}} g_2C_{\varepsilon_{3}} t^{\varepsilon_3}=t^{ \varepsilon_{1}}C_{-\varepsilon_{1}} g_1 t^{ \varepsilon_{2}} g_2 C_{\varepsilon_{3}}t^{\varepsilon_3}.$$ Now we cross out $ t^{\varepsilon_1}$ and $t^{\varepsilon_3}$ at both sides of the above equation. From the equality we obtain we deduce that there exist $d_1 \in C_{\varepsilon_1}$  and $d_{2} \in C_{\varepsilon_3}$  
such that 
$g_0 t^{ \varepsilon_{2}} g_1=d_1g_0 at^{ \varepsilon_{2}} g_1d_2$.  
By \theoc{main hnn}, there exists $c \in C_{\varepsilon_2}$ such that 
\begin{align}\label{b1}
g_0&=d_1g_0 a c\\
g_1&=\varphi^{ \varepsilon_{2}}(c^{-1})g_1 d_2.\label{b2}
\end{align}
By malnormality, separability, and Equation (\ref{b2}), $\varphi^{ \varepsilon_{2}}(c^{-1})=1$ and so $c=1$. By 
malnormality, separability and Equation (\ref{b1}), $ac=1$. Hence, $a=1$.

Now we prove $(2)$. Observe that if the two sets of the statement of the lemma are equal, then there are two possibilities:
\begin{romlist}
\item  The first (resp. second) element listed in the right set is equal to the first (resp. second) element listed on the 
second set.
\item The first (resp. second) element listed in the right set is equal to the second (resp. first) element listed on the 
second set.
\end{romlist}
The case (i) is ruled out by statement $(1)$. Hence, we can assume that (ii) holds. By \theoc{main hnn}, for each $i  \in \{1,2,3\}$, $\varepsilon_{i}=\eta_{i}$.  By arguing as in the proof of statement $(1)$, we can deduce that  there exist $c_{1}$ and $d_{1}$ in $C_{-\varepsilon_{1}}$ and  $c_{3}$ and $d_3$ in  $C_{\varepsilon_{3}}$,  such that
 $$g_1a t^{\varepsilon_2} g_2=c_1h_1 b t^{ \varepsilon_2} h_2c_{3} 
 \mbox{ and } g_1 t^{\varepsilon_2} g_2=d_{1}h_1 t^{ \varepsilon_2} h_2d_3$$

By \theoc{main hnn},   there exist a pair of elements , $x$ and $y$  in $C_{\varepsilon_{2}}$ such 
that 
\begin{align}\label{bb1}
g_0a&=c_1h_0 b x\\
g_1&=\varphi^{\varepsilon_2}(x^{-1})h_1 c_{3}.\label{bb2}\\
g_0&=d_{1}h_0  y\label{bb3}\\
g_1&=\varphi^{\varepsilon_2}(y^{-1})h_1d_{3}.\label{bb4}
\end{align}
Since $\varphi$ is an isomorphism, by malnormality and separability and Equations (\ref{bb2}) and (\ref{bb4}), $x=y$. Analogously, by 
malnormality and separability and Equations (\ref{bb1}) and (\ref{bb3}), $bx a^{-1}=y$. Therefore, $a$ and $b
$ are conjugate by $x$. Since $x \in A \cup B$,  the proof is complete.
\end{proof}

The next theorem gives necessary conditions for certain products of cyclically reduced sequences to not be 
conjugate.

\start{theo}{non-separating hnn}  Let  $G^{\ast\varphi}$ be an HNN
extension. Let $(g_0,t^{\varepsilon_1},g_1,t^{\varepsilon_2}\ldots,
g_{n-1},t^{\varepsilon_n})$ be a cyclically reduced sequence.
 and let  $i, j \in \{1,2,\ldots, n\}$. Let $a$ be an element of $C_{\varepsilon_{i}}$ and  let $b$ be an element of $C_{\varepsilon_{j}}$.
Moreover, assume that the following conditions hold:
\begin{numlist}
\item The subgroups $A$ and $B$ are malnormal in $G$.
\item The HNN extension $G^{\ast\varphi}$  is separated.
\item If $\varepsilon_i=\varepsilon_j$ then $a$ and $b$ are not conjugate by an element of $A \cup B$.
\item  Either $a \ne 1$ or $b \ne 1$.
\end{numlist}
 Then the products
$$
g_0t^{\varepsilon_1} g_1 t^{\varepsilon_2}\cdots g_{i-1}a t^{\varepsilon_{i}}g_{i}\cdots
g_{n-1}t^{\varepsilon_n} \mbox{ and } g_0 t^{\varepsilon_1}g_1 t^{\varepsilon_2}\cdots g_{j-1} b
t^{\varepsilon_{j}} g_{j}\cdots g_{n-1}t^{\varepsilon_n}
$$
are not conjugate.
\end{theo}

\begin{proof}
 Assume that the products of the hypothesis of the theorem are conjugate. Furthermore, we assume
 that $a$ is not the identity. The case of $b$ not the identity can be treated similarly.

 If $n=1$ the proof of the result is direct. Hence, we can assume $n \ge 2$.
 
The sequences
$$(g_0,t^{\varepsilon_1},g_1,t^{\varepsilon_2}\ldots,g_{i}a,t^{\varepsilon_{i+1}},
\ldots g_{n-1},t^{\varepsilon_n}) \mbox{ and }
(g_0,t^{\varepsilon_1},g_1,t^{\varepsilon_2}\ldots,g_{j}b,t^{\varepsilon_{j+1}},\ldots
g_{n-1},t^{\varepsilon_n})
$$  are cyclically reduced. By  \lemc{claim two} both sequences are  associated with the same cyclic sequence of 
double cosets. If $i=j$ the result is a consequence of \lemc{claim three}(1). Hence, we can assume that $i \ne j$. As  in the proof of \theoc{amalgamating with cyclic} by crossing out the elements in the sequences 
of double cosets with equal expressions we obtain the two sets below are equal,
$$
\{C_{\varepsilon_{i-1}} t^{\varepsilon_{i-1}} g_{i-1}a t^{\varepsilon_{i }} g_{i }  t^{\varepsilon_{i+1}} C_{-\varepsilon_{i+1}}, 
  C_{\varepsilon_{j-1}} t^{\varepsilon_{j-1}} g_{j-1}   t^{\varepsilon_{j }} g_{j }  t^{\varepsilon_{j+1}} C_{-\varepsilon_{j+1}}\}$$
$$\{C_{\varepsilon_{i-1}} t^{\varepsilon_{i-1}}g_{i-1} t^{\varepsilon_{i }} g_{i } t^{\varepsilon_{i+1}} C_{-\varepsilon_{i+1}}, 
 C_{\varepsilon_{j-1}}t^{\varepsilon_{j-1}} g_{j-1} b t^{\varepsilon_{j }}g_{j }  t^{\varepsilon_{j+1}}C_{-\varepsilon_{j+1}} \}
$$
(By \remc{epsilon},
$$C_{\varepsilon_{i-2}} t^{\varepsilon_{i-2}} g_{i-2} t^{\varepsilon_{i-1}} g_{i-1}a  t^{\varepsilon_{i}} C_{-\varepsilon_{i}}=
C_{\varepsilon_{i-2}} t^{\varepsilon_{i-2}} g_{i-2} t^{\varepsilon_{i-1}} g_{i-1}  t^{\varepsilon_{i}} C_{-\varepsilon_{i}}.
$$ Hence, double cosets like those on the left hand side of the above equation  do not appear in the list of possibly distinct double cosets).
By \lemc{claim three}, $\varepsilon_i=\varepsilon_j$ and $a$ and $b$ are conjugate by an element of $A \cup B$. 
This contradicts our hypothesis and so the proof is complete.
\end{proof}

\start{rem}{counter?} 
The following example shows that hypotheses (1) or (2) of \theoc{non-separating hnn}  are necessary. Let $G$ be the direct sum of $\Z$ and $\Z$, $\Z \oplus \Z$. Let $A = \Z \oplus \{0\}$ and let $B = \{0\} \oplus \Z$ considered as subgroups of $\Z$. Define $\map{\varphi}{A}{B}$ by $\varphi(x,0)=(0,x)$. Let $a$ be an element of $A$.   Since $\Z \oplus \Z$ is commutative,
$$
(0,1) t^{-1} (1,1) t^{-1} (1,1) t (0,-1) =t^{-1}  (2,1)   t^{-1} (0,1) t $$

Thus the sequences $t^{-1}, (1,1), t^{-1}, (0,1)(1,0), t $ and $t^{-1},  (1,1)(1,0),   t^{-1} ,(0,1), t$ have conjugate products. On the other hand, hypotheses $(3)$ and $(4)$ of \theoc{non-separating hnn} hold for these sequences.

The following example shows that hypothesis (3) of \theoc{non-separating hnn} is necessary. Let $G^{\ast}$ be an HNN extension relative to an isomorphism $\map{\varphi}{A}{B}$. Let $g \in G \setminus A$ and let $a \in A$. The sequence $(g,t,g,t)$ is cyclically reduced but the product of the sequences 
 $(ga,t,g,t)$ and $(g,t,ga,t)$ are conjugate.
 
 \end{rem}
The next auxiliary lemma will be used in the proof of \theoc{hnn inverse}. The set of congruence classes modulo $n
$ is denoted by $\mathbb{Z}/n\mathbb{Z}$.
\start{lem}{mod} Let $n$ and $k$ be integers. Let $\map{\widehat{F}}{\mathbb{Z}/n\mathbb{Z}}{\mathbb{Z}/n
\mathbb{Z}}$  be  induced by the map on the integers $F(x)=-x+k$. If $n$ is odd or $k$ is even then $\widehat{F}$ 
has a fixed point.
\end{lem}
\begin{proof}  If $k$ is even then $\frac{k}{2}$ is integer and a fixed point of $F$.  Thus $\widehat{F}$ has a fixed 
point. 

On the other hand, the map $\widehat{F}$ has a fixed point whenever the equation  $2x  \equiv -k \pmod n$ has a 
solution. If $n$ is odd this equation has a solution because $2$ has an inverse in $\mathbb{Z}/n\mathbb{Z}$. This 
completes the proof. \end{proof}

We will use \notc{C} for the statement and proof of the next result.

\start{theo}{hnn inverse}  Let  $G^{\ast\varphi}$ be an HNN
extension. Let
$(g_0t^{\varepsilon_1} g_1 t^{\varepsilon_2}\ldots g_{n-1}t^{\varepsilon_n}) $ be a cyclically reduced sequence. Let 
$i$ and $j$ be elements of  $\{1,2,\ldots, n\}$. Let  $a$ be an element of $C_{\varepsilon_{i}}$ and let $b$ an 
element of $C_{-\varepsilon_{j}}$.   Assume that for each $g \in G$, $g^{-1} $ does not belong to the set  $(A \cup 
B) g (A  \cup B)$. 
 Then the products
$$
g_0t^{\varepsilon_1} g_1 t^{\varepsilon_2}\cdots g_{i-1}a
t^{\varepsilon_{i}}g_{i}\cdots g_{n-1}t^{\varepsilon_n} \mbox{
and }$$ $$
g_{n-1}^{-1}t^{-\varepsilon_{n-1}}g_{n-2}^{-1}t^{-\varepsilon_{n-2}}\cdots
g_{j}^{-1}b t^{-\varepsilon_{j}} \cdots
t^{-\varepsilon_1}g_0^{-1}t^{-\varepsilon_n}
$$
are not conjugate.
\end{theo}
\begin{proof}
We start by giving a sketch of the proof: If the above products are conjugate then the sequence of $(\varepsilon_
{1}, \varepsilon_{2}, \dots, \varepsilon_{n})$ is a rotation of the sequence $(-\varepsilon_{n},- \varepsilon_{n-1}, 
\dots, -\varepsilon_{1})$. Since the terms of those sequences are not zero, a term of the first of the form $
\varepsilon_{h}$  cannot correspond to a term of the form  $-\varepsilon_{h}$. This gives conditions of the rotation 
and the number of terms. The same rotation also establishes a correspondence between double cosets of $g_{h}
$'s and double cosets of rotated $g_{h}^{-1}$'s. Using the fact that the sequences  $\varepsilon_{h}$'s and $g_{h}$'s are ``off'' by one, we will show that 
there exists $u$ such that $g_{u}^{-1} \in (A \cup B) g_{u} (A  \cup B)$. (in the detailed proof, we need to study 
separately the cases where the elements $a$ and $b$ appear in the equations.)

Here is the detailed proof. Consider the sequence  $(s_0,t^{\eta_1},s_1,t^{\eta_2}\ldots, s_{n-1},t^{\eta_n})$ 
defined by  $\eta_h=-\varepsilon_{n-h}$ and   
\begin{equation}\label{s}
{s_h = }\begin{cases}
g^{-1}_{-h-1} b &  \mbox{if $-h-1 \equiv j \pmod n$,}\\
g^{-1}_{-h-1}  & \mbox{otherwise.}
\end{cases}
\end{equation}
for each $h$. (Recall \notc{permutations hnn}.)

Assume that the products of the hypothesis of the theorem are conjugate.
Thus the sequences
$$(g_0,t^{\varepsilon_1},g_1,t^{\varepsilon_2}\ldots,g_{i-1}a,t^{\varepsilon_{i}}, g_{i}, 
\ldots g_{n-1},t^{\varepsilon_n}) \mbox{ and }
(s_0,t^{\eta_1},s_1,t^{\eta_2}\ldots s_{n-1},t^{\eta_n})$$ are cyclically 
reduced and have conjugate products. Let  $k$ be as in \coryc{hnn cory cory} for these two products. Hence,
\begin{equation}\label{eta}
\varepsilon_{h}=\eta_{h+k}=-\varepsilon_{n-h-k},
\end{equation}  and the following holds:
\begin{numlist}
\item if $h \not\equiv i-1 \pmod  n$ then  $g_h$ belongs to the double coset $C_{-\varepsilon_{h}}s_{h+k}  C_
{\varepsilon_{h+1}}$
\item  $g_{i-1}a$ belongs to the double coset $C_{-\varepsilon_{i-1}}s_{i+k}  C_{\varepsilon_{i}}$. \end{numlist}

By hypothesis, $a \in C_{\varepsilon_{i}}$.
Thus, for all $h$ we have
\begin{equation}\label{aaa}
g_{h} \in C_{-\varepsilon_{h}}s_{h+k}  C_{\varepsilon_{h+1}}
\end{equation}

Since for  every integer $h$,  $\varepsilon_h \ne 0$ we have that  $\varepsilon_h \ne -\varepsilon_h$. Therefore, by 
Equation (\ref{eta}), the map $\widehat{F}$ defined on the integers mod $n$ by the formula $F(h)=-h-k$ cannot 
have fixed points. By \lemc{mod}, $n$ is even and $k$ is odd. 

Then $(-k-1)$ is even. By \lemc{mod}, the map $G(h)=-h+(-k-1)$ has a fixed point. Denote this fixed point by $u$. 
Thus
\begin{equation}\label{uuu}
u+k   \equiv -u-1 \pmod n
\end{equation}
Assume that $u\not\equiv j \pmod n$. By Equations (\ref{uuu}) and (\ref{s}), $s_{u+k}=s_{-u-1}=g_{u}^{-1}$. By 
Equation \ref{aaa}, $g_{u} \in C_{-\varepsilon_{u}}g^{-1}_{u}  C_{\varepsilon_{u+1}}\subset (A \cup B) g_{j} (A  \cup 
B)$, contradicting our hypothesis. Therefore, $u\equiv j \pmod n$. In this case, by Equation \ref{aaa},  $g_{j} \in C_
{-\varepsilon_{j}}g^{-1}_{j}b  C_{\varepsilon_{j+1}}$. 

By Equations (\ref{uuu})  and (\ref{eta}), $\varepsilon_{j+1}=-\varepsilon_{-j -1-k}=-\varepsilon_{j}$. Since $b \in C_
{-\varepsilon_{j}}$ then $b \in C_{\varepsilon_{j+1}}$. Consequently, $g_{j} \in C_{-\varepsilon_{j}}g^{-1}_{j} C_
{\varepsilon_{j+1}}$. Since  $C_{-\varepsilon_{j}}g^{-1}_{j} C_{\varepsilon_{j+1}}\subset (A \cup B) g_{j} (A  \cup B)
$, this is a contradiction.

\end{proof}

\section{Oriented  non-separating simple loops}\label{non-separating case}

This section is the "separating version" of Section \ref{separating oriented}. The main purpose here consists in 
proving \theoc{bracket form hnn}, which describes the terms of the bracket of a simple non-separating conjugacy class   and an arbitrary  conjugacy class.

\begin{figure}[ht]
\begin{center}
\begin{pspicture}(14,5)

\pscurve[linecolor=gray,arrowsize=7pt]{->}(1,1)(1.8,0.8)(2.2,0.8)(3,1)
\pscurve[linestyle=dashed, linecolor=gray](1,1)(2,1.2)(3,1)
\psline(1,1)(1,4) \psline[](3,1)(3,4)
\pscurve[linecolor=gray,arrowsize=7pt]{->}(1,4)(1.8,3.8)(2.2,3.8)(3,4)
\pscurve[linecolor=gray](1,4)(2,4.2)(3,4)
\psline[linecolor=lightgray,arrowsize=5pt,linewidth=1.5pt]{<-}(2,0.8)(2,2.5)
\psline[linecolor=lightgray,linewidth=1.5pt](2,2)(2,3.8)
\psdot(2,0.77) \psdot(2,3.77) \uput[0](2,0.5){$(q,1)$}
\uput[0](1.48,3.5){(q,0)}

\pscurve[linecolor=gray,arrowsize=7pt]{->}(3.4,2.5)(4.4,2.8)(5.4,2.5)
\uput[0](4.4,3){$\psi$}

\psccurve(6,2)(8,0.2)(10,0.2)(12,0.3)(13.5,2)(10,4)

\pscurve[linecolor=gray,arrowsize=7pt]{->}(10,0.2)(9.8,1)(10,2)
\pscurve[linecolor=gray,linestyle=dashed](10,0.2)(10.2,1)(10,2)

\pscurve[linecolor=gray,arrowsize=7pt]{->}(10,4)(9.8,3)(10,2)
\pscurve[linecolor=gray,linestyle=dashed](10,4)(10.2,3)(10,2)

\rput(10,0){$\lambda_0$} \rput(10,4.25){$\lambda$}

\psecurve[linecolor=lightgray,linewidth=1.5pt](11,1)(10,2)(8,3)(7,2)(8,1)(10,2)(11,3)
\psline[linecolor=lightgray,linewidth=1.5pt,arrowsize=7pt]{<-}(7,1.8)(7,2)
\uput[0](6.7,2.5){$\tau$}

\pscurve(10,2)(9,1.79)(8,2) \pscurve(10,2)(9,2.12)(8.18,1.95)

\rput(2.6,0.6){\pscurve(10,2)(9,1.79)(8,2)
\pscurve(9.8,1.94)(9,2.12)(8.18,1.95)}
\rput(2.8,-0.6){\pscurve(10,2)(9,1.79)(8,2)
\pscurve(9.8,1.94)(9,2.12)(8.18,1.95)}

\rput(10.2,2){$p$} \psdot(10,2)
\end{pspicture}\caption{The map $\psi$ of \lemc{hnn construction}}\label{pinched}

\end{center}
\end{figure}

We will start by proving some elementary auxiliary results.

\start{lem}{hnn construction} Let $\lambda$ be  a non-separating
simple curve on $\Sigma$ and $p$ a point in $\lambda$. There exists
a map $\map{\psi}{\SI \times [0,1]}{\Sigma}$, such that:
\begin{numlist}

\item There exists a point $q \in \SI$, $\psi(q,0)=\psi(q,1)=p$.
\item $\psi$ is injective on $(\SI \times [0,1] ) \setminus \{(q,0),(q,1)\}$
\item $\psi|_{\SI \times\{0\}}=\lambda$.
\end{numlist}

\end{lem}

\begin{proof} Choose a simple  trivial  curve $\tau$ such that
$\tau \cap \lambda=\{p\}$ and the intersection of $\tau$ and
$\lambda$ is transversal. (The existence of such a curve is
guaranteed by the following argument of Poincar\'{e}: take a small
arc $\beta$ crossing $\lambda$ transversally. Since $\Sigma
\setminus (\beta \cup \lambda)$ is connected there exists a
non-trivial arc  in $\Sigma \setminus (\beta \cup \lambda)$, with no
self-intersections, joining the endpoints of $\beta$.)


\begin{figure}[ht]
\begin{center}
\begin{pspicture}(14,5)

\pscurve[linecolor=gray,arrowsize=7pt]{->}(1,1)(1.8,0.8)(2.2,0.8)(3,1)
\pscurve[linestyle=dashed, linecolor=gray](1,1)(2,1.2)(3,1)
\psline(1,1)(1,4) \psline[](3,1)(3,4)
\pscurve[linecolor=gray,arrowsize=7pt]{->}(1,4)(1.8,3.8)(2.2,3.8)(3,4)
\pscurve[linecolor=gray](1,4)(2,4.2)(3,4)
\psline[linecolor=lightgray,arrowsize=5pt,linewidth=1.5pt]{<-}(2,0.8)(2,2.5)
\psline[linecolor=lightgray,linewidth=1.5pt](2,2)(2,3.8)
\psdot(2,0.77) \psdot(2,3.77) \uput[0](2,0.5){$(q,1)$}
\uput[0](1.48,3.5){(q,0)}

\pscurve[linecolor=gray,arrowsize=7pt]{->}(3.4,2.5)(4.4,2.8)(5.4,2.5)
\uput[0](4.4,3){$\eta$}

\psccurve(6,2)(8,0.2)(10,0.2)(12,0.3)(13.5,2)(10,4)

\pscurve[linecolor=gray,arrowsize=7pt]{->}(10,4)(9.8,3)(10,2)
\pscurve[linecolor=gray,linestyle=dashed](10,4)(10.2,3)(10,2)

 \rput(10,4.25){$\lambda$}

\rput(-1.5,0){
\pscurve[linecolor=gray,arrowsize=7pt]{->}(10,4)(9.8,3)(10,2)
\pscurve[linecolor=gray,linestyle=dashed](10,4)(10.2,3)(10,2)

}

\psccurve[linecolor=lightgray,linewidth=1.5pt](8,3)(7,2)(8,1)(9,0.861)(10,1)(10.5,2)(9.8,3)
\psline[linecolor=lightgray,linewidth=1.5pt,arrowsize=7pt]{<-}(7,1.8)(7,2)
\uput[0](6.7,2.5){$\tau$}\psdot(9.8,3)

\pscurve(10,2)(9,1.79)(8,2) \pscurve(10,2)(9,2.12)(8.18,1.95)

\rput(2.5,0.6){\pscurve(10,2)(9,1.79)(8,2)
\pscurve(9.8,1.94)(9,2.12)(8.18,1.95)}

\rput(2.8,-0.6){\pscurve(10,2)(9,1.79)(8,2)
\pscurve(9.8,1.94)(9,2.12)(8.18,1.95)}

\rput(8.4,4.2){$\chi$}
\psdot*(8.3,3.1) \rput(8,3.2){$s$} \psdot*(8.34,3.46)
\rput(8.12,3.57){$s_1$} \psdot*(8.34,2.7) \rput(8,2.6){$s_2$}

\pscurve(8.12,3.57)(8.6,3.2)(8.34,2.7) \rput(8.8,3.2){$\kappa$}

\psecurve(9,3)(8.34,3.46)(6.4,2)(9,0.4)(10.74,1.53)(9.8,3)(8,2)
\psecurve(9,3)(8.34,2.7)(7.6,2)(9,1.4)(10,1.6)(9.8,3)(8,2)

\end{pspicture}\caption{The
proof of \lemc{hnn construction}}\label{lemma}

\end{center}

\end{figure}

Consider an injective map  $\map{\eta}{\SI\times [0,1]} { \Sigma}$
such that $\lambda = \eta(\SI\times \{0\})$. Let $q \in \SI$ be such
that $\eta(q,0)=p$. Denote by $C$ the image of the cylinder
$\eta(\SI\times [0,1])$. Denote by $\xi$ the boundary component of
$C$ defined by $\eta(\SI\times \{1\})$. Modify $\eta$ if necessary
so that $\tau$ intersects $\xi$ in a unique double point $s$. (See
Figure \ref{lemma}.)

Choose two points on $s_1$ and $s_2$ on $\xi$ close to $s$ and at
both sides of $s$. Choose a embedded arc in the interior of $C$,
intersecting $\tau$ exactly once,  from $s_1$ to $s_2$ and denote it
by $\kappa$. Denote by $D$ the closed half disk bounded by $\kappa$
and the subarc of $\xi$ containing $s$.

Choose two disjoint embedded arcs $\alpha_1$ and $\alpha_2$  on $\Sigma$ from $p$ to $s_1$ and  $s_2$, 
respectively, and such that $\alpha_1 \cap \tau=\alpha_2 \cap \tau=\{p\}$,  $\alpha_1 \cap C$=$\{p,s_1\}$ and $
\alpha_2 \cap C$=$\{p,s_2\}$.

 Consider the triangle $T$, with sides $\alpha_1$, the arc in $\tau$ from $s_1$ to
$s_2$ and $\alpha_2$.

Denote by $\eta_1$ the restriction of $\eta$ to $\eta^{-1}(D)$, $\map{\eta_1}{\eta^{-1}(D)}{D}$. Now, take a
homeomorphism $\map{\eta_2}{D}{D \cup T}$ such that $\eta_2|_{\kappa}$ is the identity and $\eta_{2}(s)=p$.

For each $x \in \SI\times [0,1]$, define $\map{\psi}{\SI\times
[0,1]}{\Sigma}$ as $\eta(x)$ if $x \notin D$ and as $\eta_2\eta(x)$
otherwise. This maps satisfies the required properties.

\end{proof}

Let $\psi$ be the map of \lemc{hnn
construction}. Denote by $\lambda_1$ the curve $\psi(\SI\times
\{0\})$. The  homeomorphism
$\map{\vartheta}{\lambda=\psi(\SI\times\{0\})}{\lambda_1=\psi(\SI\times
\{1\})}$ defined by $\vartheta(\psi(s,0))=\psi(s,1)$  induces an
isomorphism $\map{\varphi}{\pi_1(\lambda,p)}{\pi_1(\lambda_1,p)}$.
Denote by $\Sigma_1$ the subspace of $\Sigma$ defined by $\Sigma
\setminus \psi(\SI \times (0,1))$, and by $\tau$ the simple closed
curve induced by the restriction of $\psi$ to $\{q\} \times [0,1]$.

\start{lem}{hnn construction 2} With the above notation, the fundamental group of $\Sigma$,  $\pi_1(\Sigma,p)$, is 
isomorphic to the HNN extension of $\pi_1(\Sigma_1,p)$ relative to
$\varphi$.
Moreover,  $\tau$ is a representative of the element denoted by the stable letter $t$   and if $(g_0, t^
{\varepsilon_1},g_1,t^{\varepsilon_2},\ldots$ $g_{n-1},t^{\varepsilon_n} )$ is a cyclically reduced sequence  of the HNN extension then there exists a sequence of curves $(\gamma_0, \gamma_1,\ldots,\gamma_{n-1})$ such that for each $i \in \{0,1,\ldots,n-1\}$,
\begin{numlist}
\item The basepoint if $\gamma_i$ is $p$.
\item The curve $\gamma_i$ is a representative of $g_{i}$.
\item The inclusion $\gamma_i \subset \Sigma_1$ holds.
\item The basepoint $p$ is not self-intersection point of $\gamma_i$. In other words, $\gamma_i$ passes through
$p$ exactly once.
\end{numlist}
\end{lem}
\begin{proof} By the van Kampen's Theorem (see, for instance \cite{h}), since $\Sigma_1
\cap \tau=\{p\}$ we have that $\pi_1(\Sigma_1 \cup \tau,p)$ is the free product of
$\pi_1(\Sigma_1,p)$
 and the infinite cyclic group $\pi_1(\tau,p)$.

Denote by $D$ the disk $\psi((\SI\setminus q) \times (0,1) )$. Glue the boundary of  $D$ to the boundary of  $
\Sigma_1 \cup \tau $  as follows: attach
$\lambda$, $\tau$, $\lambda_1$ to $\SI \times \{0\}$, $q \times [0,1]$, and $\SI \times \{1\}$ respectively.
 (The reader can easily deduce the orientations.)

The relation added by attaching the disk $D$ shows that the $\pi_1(\Sigma,p)$ is isomorphic to the
HNN extension of $\pi_1(\Sigma_1,p)$ relative to $\varphi$. Notice also that $\tau$ is a
representative of $t$.

For each $i \in \{0,1,\ldots,n-1\}$,  let $\gamma_i$ be a loop in
$\Sigma_1$,  based at $p$ and representing $g_i$. By modifying these
curves by a homotopy relative to $p$ if necessary, we can assume
that each of them intersects $p$ exactly once, as desired. Then
$(4)$ follows.
\end{proof}

Recall that there is a canonical isomorphism between  free homotopy classes of curves on a surface $\Sigma$ and 
conjugacy classes of elements of $\pi_{1}(\Sigma)$. From now on, we will identify these two sets.

The following theorem gives a combinatorial description of the bracket of two oriented curves, one of them simple 
and non-separating.

\start{theo}{bracket form hnn} Let  $\lambda$ be a separating simple closed curve. Let $(g_0,t^{\varepsilon_1},g_1,
$
$t^{\varepsilon_2}\ldots,$ $g_{n-1},t^{\varepsilon_n})$ be a
cyclically reduced sequence for the HNN extension of \lemc{hnn
construction 2} determined by $\lambda$. Let $y$  be the element of $\pi_1(\Sigma,p)$
represented by $\lambda$. Then the following holds.
\begin{numlist}
\item  There exists a a representative $\eta$ of the conjugacy class of the product $g_0
t^{\varepsilon_1}g_{1}t^{\varepsilon_2} \cdots g_{n-1}$  such that $\eta$ and $\lambda$  intersect transversally at 
$p$ with multiplicity $n$.

\item  There exists $s \in \{1,-1\}$ such   the bracket is given by

\begin{align*}
s \hspace{0.1cm} [g_0 t^{\varepsilon_1}g_1 t^{\varepsilon_2}\cdots
g_{n-1}t^{\varepsilon_n},y]=&  \sum_{\textit{ $i$ :
$\varepsilon_i=1$}} g_0 t^{\varepsilon_1}g_1 t^{\varepsilon_2}\cdots
g_{i-1}\hspace{0.1cm}y \hspace{0.1cm} t^{\varepsilon_i} g_i
\cdots g_{n-1}t^{\varepsilon_n} - \\
&\sum_{ \textit{ $i$ : $\varepsilon_i=-1$} } g_0
t^{\varepsilon_1}g_1 t^{\varepsilon_2}\cdots
g_{i-1}\hspace{0.1cm}\varphi(y)\hspace{0.1cm} t^{\varepsilon_i}
g_i\cdots g_{n-1}t^{\varepsilon_n}
\end{align*}
\end{numlist}
\end{theo}

\begin{proof} Let  $(\gamma_0, \gamma_1, \ldots, \gamma_{n-1})$  denote a  sequence of  curves obtained in
\lemc{hnn construction 2} for the sequence $(g_0,t^{\varepsilon_1},g_1$,$t^{\varepsilon_2}\ldots,$
$g_{n-1}$,$t^{\varepsilon_n})$.

\begin{figure}[ht]\begin{center}
\begin{pspicture}(12,5)

\pswedge[fillstyle=solid,linecolor=white,fillcolor=lightgray](6,2.5){2.3}{-30}{30}
\pswedge[fillstyle=solid,linecolor=white,fillcolor=lightgray](6,2.5){2.3}{150}{220}
 \pscircle(6,2.5){2.3}
 \psline[arrowsize=5pt](4.2,1)(6,2.5)

 \psline[arrowsize=5pt]{->}(5.1,1.76)(5.16,1.82)

  \psline[arrowsize=5pt](6,2.5)(8,1.4)
  \psline[arrowsize=5pt]{->}(6,2.5)(4,3.6)

   \psline[arrowsize=5pt]{->}(6,2.5)(8,3.6)
   \psline[arrowsize=5pt]{<-}(7,1.96)(7.1,1.9)

\psdot*(6,2.5)\rput(6,3){$p$}

   \psline[arrowsize=5pt,linecolor=darkgray]{->}(6,2.5)(3.8,2)\rput(4,2.4){$\tau$}
\psline[arrowsize=5pt,linecolor=darkgray]{<-}(6,2.5)(8.2,2)

\rput(4,0.9){$\lambda$} \rput(3.8,3.9){$\lambda$}
\psline[arrowsize=5pt]{<-}(7,0.4)(6,2.5)

\psline[arrowsize=5pt,linecolor=darkgray,linestyle=dashed]{<-}(7,0.4)(8.2,2)
\rput(7.5,1.2){$\rho_i$} \rput(6.3,1.2){$\gamma_i$}
\rput(7.4,0.1){$q_i$} \rput(8.5,2){$r_i$} \rput(5,2.65){$T_1$}
\rput(7,2.65){$T_2$}

\end{pspicture}
\end{center}
\caption{Proof of \theoc{bracket form hnn}}\label{signs}
\end{figure}

Let $D \subset \Sigma$ be  a small disk around $p$. Observe that $D \cap \psi(\SI \times [0,1])$
consists in two "triangles" $T_1$ and $T_2$, intersecting at $p$ (see \figc{signs}.) Two of the
sides of one of these triangles are subarcs of $\lambda$. Denote this triangle by $T_1$. Suppose
that the beginning of $\tau$ is inside $T_1$ (the proof for the other possibility is analogous.)

\begin{figure}[ht]\begin{center}
\begin{pspicture}(12,5.5)

\rput(9,5.2){$\varepsilon_i=-1$} \rput(3,0){ \pscircle(6,2.5){2.3}
\pscurve[arrowsize=5pt,linecolor=gray]{->}(4.2,1)(6,2.5)(4,3.6)
\psdot*(6,2.5)\rput(6,3){$p$}
\pscurve[arrowsize=5pt,linecolor=gray]{->}(3.8,2)(6,2.5)(7,0.4)
\rput(1,1.5){$\widehat{\gamma}$}\rput(4.8,3.735){$\lambda$} }
\rput(10,1.5){$\widehat{\gamma}$}

\rput(3,5.2){$\varepsilon_i=1$} \rput(-3,0){
 \pscircle(6,2.5){2.3}
\pscurve[arrowsize=5pt,linecolor=gray]{->}(4.2,1)(6,2.5)(4,3.6)
\psdot*(6,2.5)\rput(6,3){$p$}

\pscurve[arrowsize=5pt,linecolor=gray]{<-}(3.8,2)(6,2.5)(7,0.4)

 \rput(4.8,3.735){$\lambda$} }

\end{pspicture}
\end{center}
\caption{Proof of \theoc{bracket form hnn}}\label{circle}
\end{figure}

For each $i \in \{1,2,\ldots,n\}$, if $\varepsilon_i=1$,   $\tau_i$ will denote a copy of the curve
$\tau$ and if $\varepsilon_i =-1$,   $\tau_i$ will denote a copy of the curve $\tau$ with opposite
direction.  Denote by $\gamma$ the curve $\gamma_1 \tau_1\gamma_2\tau_2\cdots \gamma_n \tau_n$. Clearly,
$\gamma$ is a representative of $g_0 t^{\varepsilon_1}g_1 t^{\varepsilon_2}\cdots
g_{n-1}t^{\varepsilon_n}$.

The intersection of $\gamma$ and $\lambda$ consists in $2n$ points, located at the beginning and
end of $\tau_i$ for each $i \in \{1,2,\ldots,n\}.$

We claim that for each $i \in \{1,2,\ldots,n\}$ if $\varepsilon_i=1$
then the intersection point of $\gamma$ and $\lambda$ located at the
end of $\tau_i$ can be removed by a small homotopy. Similarly, if
$\varepsilon_i=-1$ then the intersection  point at the beginning of
$\tau_i$ can be removed by a small homotopy.

Indeed, let $i \in \{0,1,\ldots,n\}$ be such that $\varepsilon_i=1$.
The intersection $\gamma \cap D$ contains $2n$ subarcs of $\gamma$.
Denote by $\varrho$ the subarc containing the end of $\tau_i$.
Denote by $r_i$ the  intersection of the boundary of $D$ and the end
of $\tau_i$ (see \figc{signs}). Denote by $q_i$ the intersection of
the beginning of $\gamma_{i}$ with the boundary of $D$. Choose an
arc from $r_i$ to $q_i$ which does not intersect $\lambda$ and
denote it by $\varrho_i$. In $\gamma$, replace $\varrho$ by
$\varrho_i$, (see \figc{signs}.) This proves the claim for the case
$\varepsilon_i=1$. The proof of the  case $\varepsilon_i=-1$ is
similar.

Denote by $\eta$ the curve obtained after homotoping $\gamma$ to remove the $n$ points mentioned
in the claim.

Note that $\eta$ intersects $\lambda$ at $p$ with multiplicity $n$. More
specifically, for each $i \in \{1,2,\ldots,n\}$, if $\varepsilon_i=1$
then there is an intersection at the beginning of $\tau_i$ and if
$\varepsilon_i=-1$ there in an intersection at the end of
$\tau_i$. Since $\eta$ crosses $\lambda$ these intersections are transversal. Thus, $(1)$ is proved.

Now, we will compute the bracket $[g_0 t^{\varepsilon_1}g_1 t^{\varepsilon_2}\cdots
g_{n-1}t^{\varepsilon_n},y]$ using $\eta$ and $\lambda$ as representatives. Since
$\eta$ and $\lambda$ have $n$ intersection points, we will find  $n$ terms.

Let $i \in \{1,2,\ldots,n\}$. Assume first that $\varepsilon_i=1$.
The term of the bracket corresponding to this intersection point is
obtained by inserting $\lambda$ between $\gamma_{i-1}$ and $\tau_i$.
Since the transformations we applied to $\gamma$ to obtain
$\eta$ can be now reversed, then the free homotopy class
of this term is $$g_0t^{\varepsilon_1}g_1\cdots
g_{i-1}\hspace{0.1cm}y \hspace{0.1cm} t^{\varepsilon_i} g_i\cdots
g_{n-1}t^{\varepsilon_n}.
$$

Assume now that $\varepsilon_i=-1$. The term of the bracket
corresponding to this intersection point is obtained by inserting
$\lambda$ right after $\tau_i$. This yields the element
$$
g_0t^{\varepsilon_1}g_1\cdots g_{i-1} t^{\varepsilon_i}y g_i\cdots g_{n-1}t^{\varepsilon_n}.
$$
By using the relation $t^{-1}y =\varphi(y)t^{-1}$ we see that this
element can be written as
$$g_0t^{\varepsilon_1}g_1\cdots g_{i-1}\hspace{0.1cm}\varphi(y)\hspace{0.1cm} t^{\varepsilon_i}g_i\cdots
g_{n-1}t^{\varepsilon_n}. $$

To conclude, observe that pairs of  terms corresponding to
$\varepsilon_i=1$ and $\varepsilon_i=-1$ have opposite signs because
the tangents of $\eta$ at the corresponding points point
in opposite directions, and the tangent of $\lambda$ is the same for
both terms. (see Figure \ref{circle})

\end{proof}

\section{Some results on  surface groups}\label{aux}

This section contains auxiliary results showing that certain equations do not hold in the fundamental group of the surface. These
results will be used in Sections \ref{conclusions} and \ref{unoriented loops} to prove that certain sequences satisfy the
hypothesis of Theorems \ref{amalgamating with cyclic}, \ref{amalgamating inverse}, \ref{non-separating hnn} and \ref{hnn
inverse}.

We will make use of the following well known result, (see \cite[Proposition 2.16]{LS}.)

\start{prop}{torsion} If $G$ is a free group or a free abelian group and   $g$ is an element of $G$ such  $g^n = 1$ for some non zero integer $n$ then $g$ is the identity.
\end{prop}
Let $F$ be a free group. An element of $F$ is said to be \emph{primitive} if it is a member of some basis of $F$.

\start{prop}{for power} Let $F$ be a free group and let $a$ be a primitive element of $F$   and let $A$ denote the cyclic group generated by $a$.  Then $A$ is malnormal in $F$. \end{prop}
\begin{proof} If $A$ is not malnormal, then  there exists an element $g$ in $F\setminus A$ 
 and two non-zero integers $n$ and $m$  such that  $g a^{m}g^{-1}a^{n}=1$. Denote by $\eta$    the map from $F$ to the abelianization of $F$,   $\map{\eta}{F}{F/[F,F]}$. We have that , $1=\eta(a^mga^{n}g^{-1})=(m+n)\eta(a)$. Since $a$ is primitive, $\eta(a)\ne 1$.   
 On the other hand, $F/[F,F]$ is a free abelian group. Thus, by \propc{torsion},  $m+n=0$.

Thus we found two elements of the free group $F$, $a^m$ and $g$, which commute. Hence $a^m$ and $g$ are power
of the same element $c \in F$ (see \cite{LS}, page 10,  for a
proof of this statement.) Let $k$ be an integer such that
$a=c^k$. By hypothesis, $a$  is primitive, thus $k \in \{1,-1\}$. Consequently, either $c=a$ or $c=a^{-1}.$ This implies that  $g$ is a multiple of $a$  contradicting the assumption that $g \in F\setminus A$. \end{proof}

If  $x$ is an element of the fundamental group of a surface and $x$ can be represented by a simple closed curve parallel to a boundary component then $x$ is primitive. Therefore, by \propc{for power} we have the following result.

\start{cory}{power} 
Let $\Sigma$ be an orientable surface with non-empty boundary and let $p$ be a point in $
\Sigma$. Let $a$  be an element of $\pi_{1}(\Sigma,p)$ which  can be represented by a simple closed  curve 
parallel to a boundary component of $\Sigma$. Then the cyclic group generated by   $a$ is malnormal 
in $\pi_{1}(\Sigma, {p})$.\end{cory}

\start{prop}{for power2} Let $F$ be a free group and let $a$ and $b$ be primitive elements of $F$.  If $n$ and $m$ are non-zero integers such that $a^mgb^ng^{-1}=1$ then either  $a$ and $b$ are conjugate or  $a$ and $b^{-1}$ are conjugate. \end{prop}
\begin{proof} Denote by $N$ the minimal normal subgroup of $F$ containing $b$, and by $\rho$ the quotient map $\map{\rho}{F}{F/N
}$.  Then $1 = \rho(a^mgb^ng^{-1})=\rho(a)^{m}.$
Since we added the relation $b=1$ and $b$ is primitive, $F/N$ is a free or trivial group. By \propc{torsion}, $1$ is the only element of $F/N$ of finite order. Consequently,  $\rho(a)=1$ and  $a \in N$. Then there exists an integer $k$ and $h \in F$ such that $a=h^{-1}b^k h$. Since $a$ is primitive, $k \in \{-1,1\}$,  as desired.
\end{proof}

\start{cory}{power22} Let $\Sigma_{1}$ and $\varphi$ be as in  \lemc{hnn construction 2}. If $\Sigma_{1}$ is not a 
cylinder then the HNN extension  of \lemc{hnn construction 2} is separated. 
\end{cory}
\begin{proof} Let $\lambda$ and $\lambda_1$ be as in the paragraph before \lemc{hnn construction 2}. Let $a$ and $b$ denote elements on the fundamental group of $\Sigma_1$ such that $\lambda$ and $\lambda_1$ are representatives of $a$ and $b$ respectively. Since $\Sigma_1$ is not the cylinder, $a$ and $b$ are not in the same conjugacy class. Thus, $a$ and $b$ are not conjugated. By \propc{for power2}, if $m$ and $n$ are integers then  the equation $a^mgb^ng^{-1}=1$ does not hold.  In other words, the HNN extension is separated.
\end{proof}

\start{prop}{for power3} Let $F$ be a free group and let $a$ and $g$ be elements of  $F$ such that $a$ is primitive. If $A$ denotes the cyclic group generated by $a$ then  for every $g \in F \setminus A$, the elements $g$ and $g^{-1}$ do not belong to the same double coset mod $A$.
\end{prop}

\begin{proof} Let $n$ and $m$ be two non zero integers. We will show that $g a^m g a^n \ne 1$.  Let $u = g a^{m}$. Since $g a^mga^{m}a^{-m}a^n=u^{2}a^{n-m}$. 
Thus $ u^{2}=a^{m-n}$
Since $a$ is primitive, $m-n$ is even and $u=a^{\frac{m-n}{2}}$. Hence, $g$ is a power of $a$.
\end{proof}

The following result follows straightforwardly from \propc{for power3}.

\start{cory}{power3} Let $\Sigma$ be an orientable surface with non-empty boundary and let $p$ be a point in $
\Sigma$. Let $a$ be an element of $\pi_{1}(\Sigma,p)$ such that $a$  can be represented by a simple 
closed curve freely homotopic to  a boundary component of $\Sigma$. Let $A$ denote the cyclic group generated by $a$. For every $g \in \pi_{1}(\Sigma,p) \setminus A$, $g$ and $g^{-1}$ do not belong to the same double coset mod $A$. \end{cory}

\start{prop}{power4} Let $\Sigma$ be an orientable surface with non-empty boundary which is not the cylinder. Let $p$ be a point in $
\Sigma$. Let $a$, $b$ and $g$ be elements of $\pi_{1}(\Sigma,p)$ such that $a$ and $b$ can be represented by 
simple closed curves freely homotopic to  distinct  boundary components of $\Sigma$. If $g$ is not
a multiple of $a$ nor a multiple of   $b$ then for every pair of integers  $n$ and $m$, $g a^mgb^n \ne 1$.\end{prop}

\begin{proof} Assume that  $g a^mgb^n = 1$. Notice that $m \ne 0$ and $n\ne 0$. Let $u = g a^{m}$. Since $g a^mgb^n=ga^mga^{m}a^{-m}b^n=u^{2}a^{-m}b^{n}$,
\begin{equation}\label{u square}
 u^{2}=b^{-n}a^{m}.
\end{equation}

Suppose that  $\Sigma$ has exactly two boundary components. Denote by $h$ the genus of $\Sigma
$. Since $\Sigma$ is not the cylinder, $h>1$. Then there exists a presentation of the fundamental group of $\Sigma$ such that the free generators are $a, a_
{1}, a_{2}, \ldots, a_{h},$ $b_{1}, b_{2}, \ldots, b_{h}$ and 
\begin{equation}\label{u}
b=a a_{1}b_{1}a_{1}^{-1}b_{1}^{-1} a_{2}b_{2}a_{2}^{-1}b_{2}^{-1}\cdots a_{h}b_{h}a_{h}^{-1}b_{h}^{-1}
\end{equation}

Combining Equations  (\ref{u square}) and (\ref{u}) we obtain $$u^{2}=(a a_{1}b_{1}a_{1}^{-1}b_{1}^{-1} a_{2}b_
{2}a_{2}^{-1}b_{2}^{-1}\cdots a_{h}b_{h}a_{h}^{-1}b_{h}^{-1})^{-n}a^{m}.$$
Observe that all the elements of the right hand side of the above equation are in the free generating set of the group.  We can check that both assumptions $n>0$ and $n<0$ lead to a contradiction. Since $n \ne 0$ the result is proved in this case. 

Now, assume that $\Sigma$ has three or more boundary components. Then there exists a free basis of the 
fundamental group of $\Sigma$ containing $a$ and $b$. In this case, an element of the form $a^{-m}b^{n}$ cannot 
be equal to an element of the form $u^{2}$ unless $m=0$ or $n=0$.  This concludes the proof.
\end{proof}

\section{Goldman Lie algebras of  oriented curves}\label{conclusions}

In this section we combine some of our previous results to prove \theoc{final}.

\start{defi}{terms} Let $x$ and $y$ be conjugacy classes of $\pi_1(\Sigma,p)$ such that $x$ can be represented by 
a simple loop. We associate a non-negative integer $t(x,y)$ to $x$ and $y$, called the \emph{number of terms of 
$y$ with respect to $x$} in the following way:

Firstly,  assume that $x$ has a separating representative. Let $(w_1, w_2, \ldots ,w_n)$ be cyclically reduced 
sequence for the amalgamated product
of \remc{canonical}  such that the product $w_1w_2 \cdots w_n$ is conjugate to $y$. (The existence of such a 
sequence is guaranteed by \theoc{existence}.) We define $t(x,y)=0$  if $n \le 1$  and $t(x,y)=n$ otherwise. (By 
\theoc{existence},  $t(x,y)$  is well defined if $x$ has a separating representative.)

Secondly, assume that $x$ can be represented as a non-separating closed curve.
Let $(g_0,t^{\varepsilon_1},g_1\ldots , g_{n-1},t^{\varepsilon_n})$ be a cyclically reduced sequence for the HNN 
extension defined in
\lemc{hnn construction 2} such that the product of this sequence is conjugate to $y$.  (The existence of such a 
sequence is guaranteed by \theoc{existence hnn}.) We set $t(x,y)=n$. (By \theoc{existence hnn}, $t(x,y)$ is well 
defined in this case.)
\end{defi}

Let $\alpha$ and $\beta$ be two curves that intersect transversally. The \emph{geometric intersection number of $
\alpha$ and $\beta$} is the number of times that $\alpha$ crosses $\beta$. More precisely, the geometric 
intersection number of $\alpha$ and $\beta$ is the number of pair of points $(u,v)$, where  $u$ is in the domain of 
$\alpha$, $v$ is in the domain of $\beta$, $u$ and $v$ have the same image in $\Sigma$ and the branch through 
$u$ is transversal to the branch through $v$ in the surface. Thus, the geometric intersection number of $\alpha$ 
and $\beta$ is the number of intersection points of $\alpha$ and $\beta$ counted with multiplicity.

Let $a$ and $b$ denote two free homotopy classes of curves. 
The \emph{minimal  intersection number of $a$ and $b$} denoted by 
by $i(a, b)$ is the minimal possible geometric intersection number  of pairs of curves representing $a$ and $b$.

\start{lem}{intersection} Let $x$ and $y$ be conjugacy classes of $\pi_1(\Sigma,p)$.
Assume that $x$ can be represented as a simple closed curve. Then the minimal
intersection number of $x$ and  $y$ is less than or equal to  the number of terms of $y$ with
respect to $x$. In symbols, $i(x,y) \le t(x,y)$.
\end{lem}
\begin{proof} If $x$ can be represented by a separating (resp. non-separating) curve, by \lemc{separating 
representative} (resp. \theoc{bracket form hnn}) there exist  representatives of $x$ and $y$ with exactly $t(x,y)$ 
intersection points. \end{proof}

\start{rem}{consequence}  If $x$ and $y$ are conjugacy classes of
$\pi_1(\Sigma,p)$ and $x$ can be represented as a simple closed
curve $x$, it can be proved directly that $i(x,y)=t(x,y).$ Since
this equality follows from \theoc{final}, we do not give a proof of this statement here. \end{rem}

\start{defi}{gold} Let $u$ and $v$ be conjugacy classes of
$\pi_1(\Sigma,p)$. The \emph{number of terms of the Goldman Lie bracket $[u,v]$} denoted by  $g(u,v)$ is the 
sum of the absolute values of the coefficients of the expression of
 $[u,v]$ in the basis of the vector space given by the set of conjugacy classes. \end{defi}

\start{rem}{inequality} Let $u$ and $v$ be conjugacy classes of
$\pi_1(\Sigma,p)$. Since one can compute the Lie bracket by taking
representatives of $u$ and $v$ with minimal intersection, and the
bracket may have cancellation then the number of terms of the Goldman Lie bracket $[u,v]$ is smaller or equal 
than the minimal intersection number of $u$ and $v$. In symbols, 
 $g(u,v)\le i(u,v)$.
\end{rem}

Our next result states that in the closed torus, the Goldman bracket always ``counts'' the intersection number of two free homotopy classes. 
\start{lem}{torus}Let $a$ and $b$ denote free homotopy classes of the fundamental group of  the closed oriented torus $\T$. Then $i(a,b)=g(a,b)$. Moreover, $[a,b]= \pm i(a,b)\hspace{0.2cm} a\cdot b$, where $a \cdot b$ denotes the based loop product of a representative of $a$ with a representative of $b$.\end{lem}
\begin{proof} 
Assume first  that  $a$  and  $b$ are  not a proper powers. Thus $a$  admits a simple representative $\alpha$,  and there exist a class $c$ which admits a simple representative $\delta$ and such that intersects  $\alpha$ exactly once with intersection number $+1$. Whence, $a$ and $c$ are a basis of fundamental group of $\T$.

Let $k$ and $l$ be integers such that $b=a^kc^l$ (Recall that the fundamental group of the closed torus is abelian). The universal cover of $\T$ is the euclidean plane $\R^2$. We can consider a projection map $\map{p}{\R^2}{\T}$ such that the liftings of $\alpha$ are the horizontal lines of equation $y=n$ with $n \in \Z$ and the liftings of $\delta$ are vertical lines with equation $x=n$ with $n \in \Z$. Thus there exist a representative $\beta$ of $b$ such that the liftings of $\beta$ are lines of equation $y=\frac{l}{k}x+n$ with $n \in \Z$.

One can check that the intersection of $\alpha$ and $\delta$ is exactly the projection of  $\{(0,0), (\frac{k}{l},0), (2\frac{k}{l},0), \dots,((l-1)\frac{k}{l},0)\}$. Moreover, all these points project to distinct points and the sign  of the intersection of $\alpha$ and $\delta$  at each of these points is equal to the sign of $l$. 

Thus $[a,b]=[a,a^kc^l]= l \hspace{0.2cm} a^{k+1}c^l$. Thus $g(a,b)=|l|$. Since there are representatives of $a$ and $b$ intersecting in $l$ points, $i(a,b) \le g(a,b)$. By \remc{inequality}, $i(a,b)=g(a,b)$.

Assume now  $a$ and $b$ are  classes which are proper powers. Thus there exist simple closed curves $\epsilon$ and $\gamma$ and integers $i$ and $j$ such that $\epsilon^i$ is a representative of $a$ and $\gamma^j$ is a representative of $b$. Moreover by the first part of this proof, we can assume that if  $e$ denotes the free homotopy class of $\epsilon$ and $g$ denotes the free homotopy class of $\gamma$, then $\epsilon$ and $\gamma$ intersect in exactly $i(e,g)$ points. 
Take $i$ ``parallel'' copies of $\epsilon$ very close to each other and  reconnect them so that they form a representative $\alpha$ of $a$.  Do the same with $j$ copies of $\gamma$, and denote by $\beta$ the representative of $b$ we obtain. We can perform this surgery far from the intersection points of $\epsilon$ and $\gamma$ so that the number of intersection points of $\alpha$ and $\beta$ is exactly $i \cdot j \cdot i(\epsilon, \gamma)$, whence  $(i \cdot j \cdot i(e,g)) \le i(a,b)$.
On the other hand,  the bracket is $[a,b]=\pm (i \cdot j \cdot i(e,g)) \hspace{0.2cm} (a\cdot b)$. By \remc{inequality}, $g(a,b)=i(a,b)$, as desired.

 \end{proof}

Here is a sketch of the proof of our next result (the detailed proof follows the statement): Given a free homotopy class $x$ with a simple representative and 
an arbitrary class $y$, we can write the bracket $[x,y]$ as in \theoc{bracket form} or as in \theoc{bracket form hnn}. 
The conjugacy classes of the terms of this bracket have representatives as the ones studied in  \theoc
{amalgamating with cyclic} or in \theoc{non-separating hnn}. By \coryc{power} (resp. \coryc{power3}) the free 
product on amalgamation (resp. the HNN extension) we are considering satisfies the hypothesis of 
\theoc{amalgamating with cyclic} (resp. \theoc{non-separating hnn}.) This implies that the terms of the bracket   do not cancel.

\start{theo}{final} Let $x$ be a free homotopy class that can be represented by a simple closed curve
on $\Sigma$ and let $y$ be any free homotopy  class. Then the following non-negative integers are equal.
\begin{numlist}
\item The minimal number of intersection points of $x$ and $y$.
\item The number of terms in the Goldman Lie bracket $[x,y]$, counted with multiplicity.
\item The number of terms in the reduced sequence of $y$ determined by $x$.
\end{numlist}
In symbols, $i(x,y)=g(x,y)=t(x,y).$
\end{theo}

\begin{proof} By Remarks  \ref{intersection} and \refc{inequality},  $g(x,y)\le i(x,y) \le t(x,y).$ Hence it is enough to 
prove that $t(x,y) \le g(x,y)$. By \lemc{torus}, we can assume that $\Sigma$ is not the torus. 
 
We first prove that $t(x,y) \le g(x,y)$ when  $x$ can be represented by a separating simple closed curve $\chi$. By 
\theoc{existence}, there exists a cyclically reduced sequence $(w_1,w_2,\ldots, w_n)$ for the free product of 
amalgamation determined by $\chi$ in \remc{canonical} such that  the product $w_{1}w_{2}\cdots w_{n}$
is a representative of  $y$. If $n=0$ or $n=1$, then $t(x,y)=0$ and the result holds. If
$n>1$, by \theoc{bracket form}, there exists $s \in \{1,-1\} $ such
that 
\begin{equation}\label{ij} s \hspace{0.1cm} [x,y]= \sum_{i=1}^n (-1)^{i}w_1w_2\cdots w_i
xw_{i+1}\cdots w_n.\end{equation} If $i$ and $j$ are such that there is cancellation between the $i$-th term and the 
$j$-th term of the right hand side of the Equation~(\ref{ij}), then $(-1)^i=-(-1)^j$. Consequently, $i$ and $j$ have different parity.

We will work at the basepoint indicated above and will abuse notation by pretending $x$ and $y$ are elements of the fundamental group of $\Sigma$ (see  \notc{declaration}.) By \coryc{power}, the cyclic group generated by $x$ is malnormal in $ \pi_1(\Sigma_1,p)$ and is malnormal in  $ 
\pi_1(\Sigma_2,p)$, where $\Sigma_{1}$ and $\Sigma_{2}$ are as in \remc{canonical}.  Thus, the hypothesis of  
\theoc{amalgamating with cyclic} hold for this free product with amalgamation.  Hence, by \theoc{amalgamating 
with
cyclic}, there is no cancellation in the sum of the right hand side of Equation~(\ref{ij}). Consequently, 
 $t(x,y)=g(x,y)$.

Now we prove the result when  $x$ can be represented by a non-separating simple closed curve $\lambda$. 
Consider  the HNN extension of \lemc{hnn construction 2} determined by $\lambda$. 
By \theoc{existence hnn}  there exists a cyclically reduced sequence 
$(g_0,t^{\varepsilon_1},g_1,t^{\varepsilon_2}, \ldots,$ $ g_{n-1},t^{\varepsilon_n})$ such that the product
$g_0t^{\varepsilon_1}g_1t^{\varepsilon_2}\cdots
g_{n-1}t^{\varepsilon_n}$ is a representative of  $y$.   By \theoc{bracket
form hnn}, there exists $s \in \{1,-1\}$ such that $s [x,y]$ equals
\begin{equation}\nonumber
 \sum_{\textit{ $i$ :
$\varepsilon_i=1$}} g_0 t^{\varepsilon_1}g_1 t^{\varepsilon_2}\cdots
g_{i-1}xt^{\varepsilon_i} g_i \cdots g_{n-1}t^{\varepsilon_n} - 
\sum_{ \textit{ $i$ : $\varepsilon_i=-1$} } g_0 t^{\varepsilon_1}g_1 t^{\varepsilon_2}\cdots
g_{i-1}\varphi(x)t^{\varepsilon_i}  \cdots g_{n-1}t^{\varepsilon_n}.
\end{equation}
If $t(x,y) >  g(x,y)$ then there is cancellation in the above sum. Therefore, there exist two integers $h$ and $k$ 
such that $\varepsilon_h=1$, $\varepsilon_k=-1$ and
the products 
$$ g_0 t^{\varepsilon_1}g_1 t^{\varepsilon_2}\cdots g_{h-1}x t^{\varepsilon_h} t^{\varepsilon_{h+1}}\cdots g_
{n-1}t^{\varepsilon_n} \mbox{ and } 
  g_0 t^{\varepsilon_1}g_1 t^{\varepsilon_2}\cdots
g_{k-1} \varphi(x)t^{\varepsilon_k}\cdots g_{n-1}t^{\varepsilon_n}
$$
are conjugate. We use now the notations of   \lemc{hnn construction 2} and the paragraph before  \lemc{hnn 
construction 2}.
By \coryc{power}, $\pi_{1}(\lambda,p)$ and  $\pi_{1}(\lambda_{1},p)$  are malnormal in $\pi_{1}(\Sigma_{1},p)$. 
Since we are assuming that $\Sigma$ is not a torus, then $\Sigma_{1}$ is not a cylinder. Then by \coryc{power22} 
the HNN extension is separated.  
By  \theoc{non-separating hnn} the two products above are not conjugate, a contradiction. Thus the proof is 
complete. \end{proof}

\begin{cory} If $x$ and  $y$ are conjugacy classes of curves that can be represented by simple closed curves then 
$t(x,y)=t(y,x)$.
\end{cory}

\begin{rem} Let $x$ be a free homotopy class that can be represented by a simple closed curve
on $\Sigma$ and let $y$ be any free homotopy  class. 
\end{rem}
\section{Goldman Lie algebras of unoriented curves}\label{unoriented loops}

Recall that Goldman \cite{g}  defined a Lie algebra of unoriented loops as follows. Denote by
$\pi^\ast$ the set of conjugacy classes of $\pi_1(\Sigma,p)$. For each $x \in \pi^\ast$,
denote by $\overline{x}$ the conjugacy class of a representative of $x$ with opposite orientation. Set
$\widehat{x}=x + \overline{x}$ and $\widehat{\pi}=\{\overline{x}+x, x \in \pi^\ast\}.$ The map $\widehat{\cdot}$ is 
extended linearly to the vector space of linear combinations of elements of $\pi^\ast$.
Denote by $V$ the real vector space generated by $\widehat{\pi}$, that is, the image of the map $\widehat{\cdot}$.  
For each pair of elements of $\pi^\ast$, $x$ and $y$, define the unoriented bracket
$$[\widehat{x},\widehat{y}]=\Big([x,y]+[\overline{x},\overline{y}]\Big)+\Big([x,\overline{y}]+[\overline{x},y]\Big)=
\widehat{[x,y]}+\widehat{[x,\overline{y}]}.
$$
We denote the bracket of oriented curves and the bracket of unoriented 
curves by the same symbol, $[\hspace{0.1cm},\hspace{0.1cm}]$.

An \emph{unoriented term} of the bracket $[a + \overline{a},b+\overline{b}]=[\widehat{a},\widehat{b}]$ of a pair of 
unoriented curves $\widehat{a}=a + \overline{a}$ and $\widehat{b}=b+\overline{b}$ is a term  of the form $c(z +
\overline{z})$, where $c$ is an integer and $z$ is a conjugacy class, that is, an elements of the basis of $V$ 
multiplied by an integer coefficient. 

Let $u(x,y)$ denote the number of unoriented terms of the bracket
$[\widehat{x},\widehat{y}]$ considered as a bilinear map on $V$,
counted with multiplicity. This is the sum of the absolute value of the coefficients of the expression of $[\widehat{x},
\widehat{y}]$ in the basis $\widehat{\pi}$.

\start{ex}{torus unoriented} With the notations of \lemc{torus}, 
$$[\widehat{a},\widehat{b}]=[\widehat{a},\widehat{a^kc^l}]= \pm i(a,b)(\widehat{a^{k+1}c^l}-\widehat{a^{k-1}c^{l}}).$$ 
Using the fact that every free homotopy class of curves in the torus admits a multiple of a simple closed curve as representative, we can show that for every pair of free homotopy classes $x$ and $y$, 
 $$u(x,y)=2 \cdot i(x,y) = 2 \cdot 
g(x,y) = 2 \cdot t(x,y).$$
\end{ex}

The strategy of the proof of our next theorem is similar to that of \theoc{final}, namely, we write the terms of the 
bracket in a certain form (using \theoc{bracket form} or \theoc{bracket form hnn}). By the results on Section~\ref
{aux}, we can apply  Theorems \ref{amalgamating inverse} and \ref{hnn inverse} to show that the pairs of  
conjugacy classes of these sums which have  different signs are distinct by Theorems \ref{amalgamating inverse} 
and \ref{hnn inverse}.

\start{theo}{unoriented}  Let $x$ and $y$ be conjugacy classes of
$\pi_1(\Sigma,p)$ such that $x$ can be  represented by a simple
closed curve. Then the following non-negative integers are equal
\begin{numlist}
\item The number of unoriented terms of the bracket $[\widehat{x},\widehat{y}]$, $u(x,y)$.
\item  Twice the minimal number of intersection points of $x$ and $y$, $2\cdot  i(x,y)$.
\item Twice the number of terms of the Goldman bracket on oriented curves, counted with multiplicity, $2\cdot g
(x,y).$
\item Twice the  number of terms of the sequence of $y$ with respect to $x$, $2\cdot t(x,y).$
\end{numlist}
In symbols,
$$
u(x,y)=2\cdot i(x,y) = 2 \cdot g(x,y) = 2 \cdot t(x,y).
$$
\end{theo} 
\begin{proof}
By \remc{inequality}, $u(x,y) \le 2\cdot  i(x,y) = 2\cdot  t(x,y)$.  By \exc{torus unoriented}, we can assume that $\Sigma$ is not the torus.
By \theoc{final}, it is enough to prove that $u(x,y)=2
\cdot t(x,y)$. Assume that $u(x,y)<2\cdot t(x,y)$.

By definition the bracket $[\widehat{x},\widehat{y}]$ is an
algebraic sum of terms of the form $\widehat{z}=z+\overline{z}$,
where $z$ is a conjugacy class of curves and $z$ and $\overline{z}$ are terms of one
of the four following brackets: $[x,y], [x,\overline{y}], [\overline{x},y]$ and $[\overline{x},\overline{y}]$.
 By
\theoc{final}, the number of terms of each of the four brackets
above is $t(x,y)$. Since $u(x,y)<2\cdot t(x,y)$ there has to be one term belonging to one of above four 
brackets that cancels with a term of another of those brackets. Denote one of these terms that cancel by $t_1$ and 
the other by  $t_2 $. By \theoc{final}, there is no  inner cancellation in any of the above four brackets. 
Consequently, if $t_1$ is a term of $[u,v]$, where  $u \in \{x, \overline{x}\}$ and $v \in \{y, \overline{y}\}$ then $t_{2}
$ is a term of  of one of the following brackets: $[u,\overline{v}]$,  $[\overline{u},v]$, or $[\overline{u},\overline{v}]$. 
Hence it suffices to
assume that $t_1$ is a term of $[x,y]$ and to analyze each of the
following three cases.
\begin{numlist}
\item  $t_2$ is a term of $[\overline{x},y]$
\item $t_2$ is a term of $[x,\overline{y}]$
\item $t_2$ is a term of $[\overline{x},\overline{y}]$
\end{numlist}

Assume first that $x$ can be represented by a separating curve $\chi$. By
\theoc{existence}, there exists a cyclically reduced sequence
$(w_1,w_2,\ldots w_n)$ in the amalgamating product of \remc{canonical} determined by $\chi$ 
such that the product $w_1w_2\cdots w_n$ is a representative of $y$.
By \theoc{bracket form},  then there exists $s \in \{1,-1\}$ and $i,
j\in \{1,2,\ldots,n\}$ such that $t_1=s (-1)^iw_1w_2\cdots w_i x
w_{i+1}\cdots w_n$ and one of the following holds.

\begin{numlist}
\item $t_2=s (-1)^{j+1}w_1w_2\cdots w_j x^{-1}w_{j+1}\cdots w_n.$
\item $t_2=s (-1)^{j+1} w_n^{-1} w_{n-1}^{-1}\cdots w_{j+1}^{-1}x w_j^{-1}\cdots
w_1^{-1}$.
\item$t_2=s (-1)^{j} w_n^{-1} w_{n-1}^{-1}\cdots w_{j+1}^{-1}x^{-1}w_j^{-1}\cdots
w_1^{-1}$.
\end{numlist}
(Note that when we change direction of  one of the elements of the bracket, $x$ or
$y$, there is a factor $(-1)$ because one of the tangent vectors at the intersection point has an opposite direction.
Also, changing direction of $x$ and $y$ does not change signs.)

Let us study first case $(1)$.   Clearly, if $t_1$ and $t_2$ cancel
then $(-1)^i=-(-1)^{j+1}$ and the products $w_1w_2\cdots
w_i x w_{i+1}\cdots w_n$ and $w_1w_2\cdots w_j x^{-1}w_{j+1}\cdots
w_n$ are conjugate. Therefore, $i$ and $j$ have equal parities. By \coryc{power}, the subgroup generated by $x$ is malnormal in the two amalgamated 
groups of the amalgamated product of \remc{canonical}. (We are again treating $x$ as an element of the 
fundamental group of the surface.) Hence we can apply  \theoc{amalgamating with
cyclic}, with $a=x$ and $b=\overline{x}$ to   show that $w_1w_2\cdots
w_i x w_{i+1}\cdots w_n$ and $w_1w_2\cdots w_j x^{-1}w_{j+1}\cdots
w_n$ are not conjugate, a contradiction.

Similarly, by \theoc{amalgamating inverse} and \coryc{power3}, Cases $(2)$ and $(3)$ are not possible.

Now assume that $x$ has a non-separating representative. By \theoc{existence hnn} there
exist a cyclically reduced sequence $(g_0, t^{\varepsilon_1}, g_1,
\ldots g_{n-1},t^{\varepsilon_{n}})$ whose product is an
element of $y$. By \theoc{bracket form hnn} there exist and integer $i$ such that the  term $t_{1}$ has the form $s
{\varepsilon_{i}} \cdot g_0 t^{\varepsilon_1} g_1 \cdots g_i u t^{\varepsilon_i}\cdots g_{n-1}t^{\varepsilon_{n}}$ 
where $u=x$ if $\varepsilon_{i}=1$ and $u=\varphi(x)$ if $\varepsilon_{i}=-1$.  By \theoc{bracket form hnn} there 
exist an integer $j$ such that the  term $t_{2}$ has one of the following forms.
\begin{numlist}
\item $t_2=-s\varepsilon_{j}\cdot g_0
t^{\varepsilon_1} g_1 \cdots g_j v t^{\varepsilon_j}\cdots
g_{n-1}t^{\varepsilon_{n}}$ where $v=x^{-1}$ if $\varepsilon_{j}=1$ and $v=\varphi(x^{-1})$ if $\varepsilon_{j}=-1$.

\item $t_2=s \varepsilon_{j} \cdot g_{n-1}^{-1}t^{-\varepsilon_{n-1}}g_{n-2}^{-1}t^{-\varepsilon_{n-2}}\cdots g_{j}^
{-1}v t^{-\varepsilon_{j}} \cdots  g_1^{-1} t^{-\varepsilon_1}g_0^{-1}t^{-\varepsilon_n}$, 
where $v=x$ if $\varepsilon_{j}=-1$ and $v=\varphi(x)$ if $\varepsilon_{j}=1$.

\item $t_2=-s \varepsilon_{j} \cdot g_{n-1}^{-1}t^{-\varepsilon_{n-1}}g_{n-2}^{-1}t^{-\varepsilon_{n-2}}\cdots
g_{j}^{-1}v^{-1}t^{-\varepsilon_{j}} \cdots  g_1^{-1}
t^{-\varepsilon_1}g_0^{-1}t^{-\varepsilon_n}$ where $v=x$ if $\varepsilon_{j}=-1$ and $v=\varphi(x)$ if $
\varepsilon_{j}=1$.

\end{numlist}

The argument continues similarly to that  of the separating case: By \coryc{power} 
and \coryc{power22} the HNN is separated and the subgroups we are considering are malnormal. 
By  \theoc{non-separating hnn}, $t_{2}$ cannot have  the form  described in Case (1). (We need to use the fact that a non-trivial 
element and its inverse are not conjugate in an infinite cyclic group). Cases $(2)$ and
$(3)$ are ruled out by \theoc{hnn inverse} and \coryc{power3} and \propc{power4}. 
\end{proof}

Let $n$ be a positive integer and let $x$ be a free homotopy class with representative $\chi$. Denote by $x^n$ the conjugacy class of the curve that wraps  $n$ times around $\chi$. We can extend \theoc{unoriented} to the case of multiple curves using the same type of arguments.

\start{theo}{multiple}  Let $n$ be a positive  integer and let $x$ and $y$ be conjugacy classes of
$\pi_1(\Sigma,p)$ such that $x$ can be  represented by a simple
closed curve $\chi$.   Then the following equalities hold.
$$
u(x^n,y)=n \cdot u(x,y) = 2 i(x^n,y) = 2\cdot n \cdot i(x,y) =  2\cdot g(x^n,y)= 2\cdot n \cdot g(x,y)=2\cdot n \cdot t(x,y).
$$
\end{theo} 

The next lemma is well known but we did not find an explicit proof in the literature. A stronger version of this result (namely, $i(x,x)$ equals twice the minimal number of  self-intersection points of $x$) is proven in \cite{CL} for the case of surfaces with boundary. 

\start{lem}{simple rep} If $x$ is a homotopy class which does not admit simple representatives then the minimal intersection number of $x$ and $x$ is not zero. In symbols, $i(x,x) \ne 0$.
\end{lem}
\begin{proof}  Let $\alpha$ and $\beta$ be two transversal representatives of $x$.  Let $\chi$ be a representative of $x$ with minimal number of self-intersection points and let $P$ be a self-intersection point of $\chi$. Let $\map{p}{H}{\Sigma}$ be the universal cover of $\Sigma$. 
Consider  two distinct lifts of $\chi$, $\widetilde{\chi}_1$ and $\widetilde{\chi}_2$ which intersect at a point $Q$ of $H$ such that $p(Q)=P$. 
Consider a lift of $\alpha$, $\widetilde{\alpha}$ such that the endpoints of $\widetilde{\alpha}$ and $\chi_1$ coincide.  Analogously, consider a lift of $\beta$, $\widetilde{\beta}$ such that the endpoints of $\widetilde{\beta}$ and $\chi_2$ coincide. Since $\chi$ intersects in a minimal number of points, the endpoints of the lifts $\chi_1$ and $\chi_2$ are linked.  Thus, the endpoints of $\alpha_1$ and $\widetilde{\beta}$ are linked. Hence, $\widetilde{\alpha}$ and $\widetilde{\beta}$ intersect in $H$. Consequently, $\alpha$ and $\beta$ intersect in $\Sigma$. 
\end{proof}

Our next result is a global a characterization of free homotopy classes with simple representatives in terms of the Goldman 
Lie bracket. 

\start{cory}{sefl} Let $x$ be a free homotopy class. Then $x$ has 
a multiple of a simple representative if and only if for every free homotopy class $y$ the number of terms of the bracket of $x$ 
and $y$ is equal to the minimal intersection number of $x$ and $y$. In symbols, $x$ has a multiple of a  simple representative if and only 
if $g(x,y)=i(x,y)$ for every free homotopy class $y$.
\end{cory}
\begin{proof} If $x$ has a simple representative and $y$ is an arbitrary free homotopy class then $g(x,y)=i(x,y)$ by 
\theoc{final}. If $x$ does not have a simple representative then by \lemc{simple rep}, $i(x,x) \ge 1$. 
On the other hand, by the antisymmetry of the Goldman Lie bracket we have $[x,x]=0$. Thus, $g(x,x)=0$  and the proof of the corollary is complete.
\end{proof}

\section{Examples}\label{examples}

The assumption in Theorems \ref{final} and \ref{unoriented} that one of the curves is simple cannot be dropped. Goldman \cite{g} gave the
following example, (attributed to Peter Scott): For any conjugacy class $a$, the Lie bracket
$[a, a]=0$. On the other hand, if $a$ cannot be represented by a multiple of a  simple curve, then any
two representatives   of  $a$ cannot be disjoint.

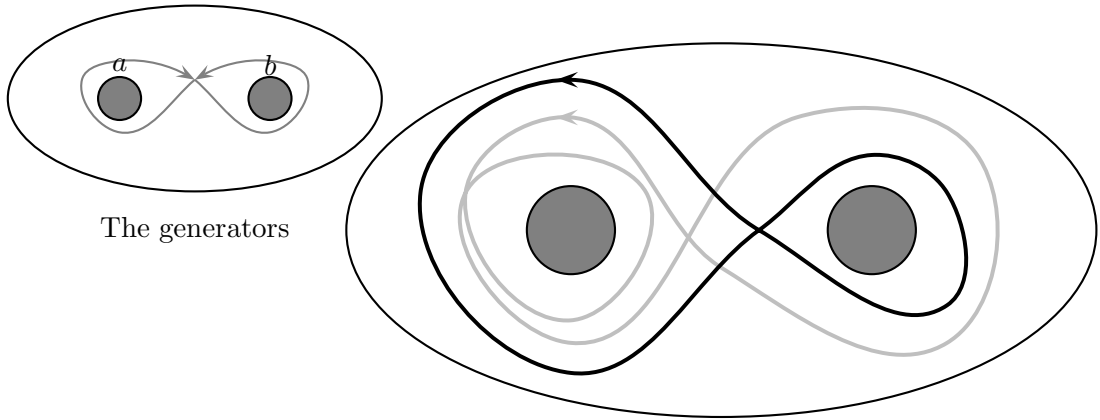
\begin{figure}[ht]\label{example}
\begin{center}
\begin{pspicture}(12,5)

\rput(-1,3){
\psset{xunit=0.5,yunit=0.5}

\psellipse[fillstyle=solid](6,2.5)(5,2.5)
\pscircle[fillstyle=solid,fillcolor=gray](4,2.5){0.3}
\pscircle[fillstyle=solid,fillcolor=gray](8,2.5){0.3}

\psecurve[arrowsize=5pt,linecolor=gray]{->}(7,4)(6,3)(4,1.6)(3,3)(4,3.5)(6,3)(7,2)

\psecurve[arrowsize=5pt,linecolor=gray]{->}(5,4)(6,3)(8,1.6)(9,3)(8,3.5)(6,3)(5,2)
\rput(4,3.4){$a$} \rput(8,3.4){$b$}
}
\rput(2,2.5){The generators}

\rput(3,0){
\psellipse[fillstyle=solid](6,2.5)(5,2.5)
\pscircle[fillstyle=solid,fillcolor=gray](4,2.5){0.6}
\pscircle[fillstyle=solid,fillcolor=gray](8,2.5){0.6}

\psline[arrowsize=6pt]{<-}(3.8,4.5)(4,4.5)

\psline[arrowsize=6pt,linecolor=lightgray]{<-}(3.8,4)(3.94,4)
\psccurve[showpoints=false,linecolor=lightgray,linewidth=1.5pt](6,2)(4,4)(2.6,3)(4,1.3)(5,3)(2.6,3)(4,1)(7,4)(9.4,3.5)
(9,1)

\psccurve[showpoints=false,linewidth=1.5pt](6.5,2.5)(4,4.5)(2,3)(4,0.6)(6.5,2.5)(8,3.5)(9,3)(9,1.5)

}
\end{pspicture}

\end{center}
\caption{An example of a pair of non-disjoint curves with bracket
$0$}
\end{figure}

Here is a family of examples:

\start{ex}{non-simple}  Consider the conjugacy classes of the curves
$aab$ and $ab$ in the pair of pants (see \figc{example}.)  The term
of the bracket correspond to $p_1$ is the conjugacy class of $aabba$
and the term of the bracket corresponding to $p_2$ is $baa ab$.The
conjugacy classes of both terms are the same, and the signs are
opposite. Then the Goldman bracket of these conjugacy classes is
zero. Nevertheless, the minimal intersection number is two (see Figure \ref{example}.)

More generally, for every pair of positive integers $n$ and $m$, the curves $a^nb$ and $a^m b$ have
minimal intersection $2\min(m,n)$. Nevertheless, the bracket of these pairs is zero. In
symbols,  $[a^nb,a^m b]=0$. (The intersection number as well as the Goldman Lie bracket can be
computed using results in \cite{chas}.)
\end{ex}

\section{Application: Factorization of Thurston's map}\label{applications 1}

Denote by $C(\Sigma)$ the set of all conjugacy classes of curves on
a surface $\Sigma$ which admit a simple representative.
Denote by
$W$ the vector space with basis all conjugacy classes of the
fundamental group of $\Sigma$. Consider the map
$\map{\phi}{C(\Sigma)}{W^{C(\Sigma)}}$, defined by
$\phi(a)(b)=[a,b]$
For each (reduced) linear combination $c$ of elements of the vector
space $W$, define a map  $\map{\mathrm{abs}}{W}{\Z_{\ge 0}} $, where $\mathrm{abs}(c)$  is the sum of the 
absolute value of the
coefficients of $c$.

By \theoc{final} the composition $\map{\mathrm{abs} \circ
\phi}{C(\Sigma)}{\Z_{\ge 0}^{C(\Sigma)}}$ is the map defined by Thurston in
\cite{T} which is used to define his mapping class group invariant compactification of Teichm\"{u}ller space.

\section{Application: Decompositions of the vector space generated by conjugacy classes}\label{applications 2}

For each  $w$ in $W$, the vector space  of free homotopy classes of curves, the \emph{adjoint map 
determined by $w$}, denoted by  $\mathrm{ad}_w$ is defined for each $y \in W$ by $\mathrm{ad}_w(y)=[y,w]$. 

Let $a$ denote the conjugacy class of a closed curve on a surface $\Sigma$ and let $n$ be a non-negative integer.
Denote by
$W_n(a)$ the subspace of $W$  generated by  the set of conjugacy classes of
oriented curves with minimal number of intersection points with $a$ equal to $n$.

In this section we prove that if $a$ is the conjugacy class of a simple closed curve  then $W_n(a)$ is invariant under $\mathrm{ad}_{a}$. Moreover, we give a further decomposition of $W_n(a)$ into subspaces invariant under $\mathrm{ad}_{a}$. 

\subsection{The separating case}\label{double cosets separating}

Let $\chi$ be a separating simple closed curve. Let $G\ast_{C}H$ be the amalgamated free product defined by $\chi$ in \remc
{canonical}. Let $(w_1, w_2, \ldots, w_n)$ be a cyclically reduced sequence for this amalgamated free product. 
Denote by $W(w_1, w_2, \ldots, w_n)$ the subspace  generated by the conjugacy classes which have  
representatives   of the form $v_{1}v_{2}\cdots v_{n}$ where for each $i$ in  $\{1,2,\dots,n\}$, $v_{i}$ is an element of 
the double coset $Cw_{i}C$. 

The following result is an immediate consequence of  \theoc{bracket form} and the definition of the subspaces $W(w_1, 
w_2, \dots, w_n)$.

\start{prop}{invariant} With the notations above, the subspace $W(w_1, w_2, \cdots, 
w_n)$ is invariant under  $\mathrm{ad}_{x}$,  the adjoint map determined by $x$. 
\end{prop}

\start{prop}{filtration mini} With the above notations,  the subspace 
$W_{n}(x)$ is the disjoint union of $W(w_{1}, w_{2}, \ldots, w_{n})$, where $(w_{1},w_{2},\ldots, w_{n})$ runs over all cyclically reduced sequences of $n$ terms of the free product with amalgamation determined by $x$. Therefore, $W_n(x)$ is invariant under $\mathrm{ad}_x$.
\end{prop}
\begin{proof} 
It is enough to prove the 
result for every element of the basis of $W_n(x)$.
Let $y \in W_n(x)$ be a conjugacy class whose minimal intersection number with $x$ is $n$.   By \theoc{existence} 
and \theoc{final}, there exists a cyclically reduced sequence $(w_1, w_2, \ldots, w_n)$ with $n$ terms
 for the amalgamating product of \remc{canonical} determined by a representative of $x$ with product in $y$.

By Theorem \ref{bracket form}, the conjugacy classes associated with each of the  terms of the bracket $[y,x]$ 
have the form
$w_1w_2 \cdots w_i x w_{i+1}\cdots w_n$ for some $i \in \{1,2,\ldots, n\}$.  The sequence $(w_1, w_2, \ldots w_ix, 
w_{i+1},\ldots , w_n)$ is cyclically reduced for every $i \in \{1,2,\ldots,n\}$. By \theoc{final}, the minimal intersection 
number of the product $w_1w_2 \cdots w_i x w_{i+1}\cdots w_n$ and $a$ is $n$. Hence each term of the Lie 
bracket $[y,x]$ is in $W_n(x)$.
\end{proof}

\start{rem}{dehn local}  It is not hard to see that the subspaces $W(w_{1},w_{2},\dots,w_{n})$ are also invariant 
under the map induced by the Dehn twist around $x$. Indeed, using \lemc{separating representative}, it is not hard 
to see that the Dehn twist around $x$ of the conjugacy class of $w_1w_2 \cdots w_n$ can be represented by 
$$w_1xw_2 \overline{x} w_{3} x w_{4}\overline{x} \cdots w_n \overline{x}\mbox{ or } w_1\overline{x}w_2xw_{3}\overline{x}\cdots w_n x$$ where the choice 
between two to conjugacy classes above is determined by the orientation of the surface, the orientation of $x$ and the subgroup which  $w_{1}$ belongs.
\end{rem}

\start{ques}{invariant subspace}  Denote by $X$ the cyclic group of automorphisms of $W$ generated by $\mathrm{ad}_x$. Let $(w_{1},w_{2},\ldots,w_{n})$ be a cyclically reduced sequence. It is not hard to see that the subspace $W(w_1, w_2, \ldots, w_n)$ is invariant under $X$. Is $W(w_1, w_2, \ldots, w_n)$ the minimal (with respect to inclusion) subspace of $W$ containing  the conjugacy class of $w_1w_2, \cdots, w_n$ and invariant under $X$.
\end{ques}

\subsection{The non-separating case}\label{double cosets non separating}

We now develop the analog of Subsection \ref{double cosets separating} for the non-separating case. In the separating case, the terms of  cyclic sequences of double cosets belong to alternating subgroups. In the non-separating case, there is a further piece of information, the sequence of $\varepsilon$'s. Thus, the arguments, although essentially the same, have some extra technical complications. Also, there are more invariant subspaces surfacing.

Let $\gamma$ denote a separating curve. Let $G^{\ast \varphi}$ be the HNN extension constructed  in \lemc{hnn construction 2}.
with $\gamma$. Let  $(g_0,t^{\varepsilon_1},g_1,t^{\varepsilon_2}, \ldots, g_{n-1},t^{\varepsilon_n})$ denote a cyclically reduced sequence. Denote by
 $W(g_0,{\varepsilon_1},g_1,{\varepsilon_2}, \ldots, g_{n-1},{\varepsilon_n})$  the subspace of $W$ generated by the conjugacy classes which have representatives of the form   $h_0t^{\varepsilon_1}h_1t^{\varepsilon_2}\cdots h_{n-1}t^{\varepsilon_n}$
where for each $i \in \{1,2,\dots,n\}$, $h_{i}$ is an element of 
the double coset $C_{-\varepsilon_{i}}w_{i}C_{\varepsilon_{i+1}}$. Let $({\varepsilon_1},{\varepsilon_2}, \ldots,{\varepsilon_n})$ be a sequence of integers such that for each $i \in \{1,2,\ldots,n\}$, $\varepsilon_{i} \in \{-1,1\}$. Denote by
 $W({\varepsilon_1},{\varepsilon_2}, \ldots,{\varepsilon_n})$  the subspace of $W$ generated by the conjugacy classes which have representatives of the form   $h_0t^{\varepsilon_1}h_1t^{\varepsilon_2}\cdots h_{n-1}t^{\varepsilon_m}$
where $(h_0,{\varepsilon_1},h_1,{\varepsilon_2},\cdots, h_{n-1},{\varepsilon_n})$ is a cyclically reduced sequence.

The following result is a direct consequence of  \theoc{bracket form hnn} and the definition of the subspaces $W(g_0,{\varepsilon_1},g_1,{\varepsilon_2}, \ldots, g_{n-1},{\varepsilon_n})$ and  $W({\varepsilon_1},{\varepsilon_2}, \ldots,{\varepsilon_n})$.

\start{prop}{invariant hnn} With the above notations, the subspaces  $W({\varepsilon_1},{\varepsilon_2}, \ldots,{\varepsilon_n})$ and  $W(g_0,{\varepsilon_1},g_1,{\varepsilon_2}, \ldots, g_{n-1},{\varepsilon_n})$ are invariant under  $\mathrm{ad}_{y}$,  the adjoint map determined by $y$. Moreover,  $W({\varepsilon_1},{\varepsilon_2}, \ldots,{\varepsilon_n})$ is a disjoint union of subspaces of the form
 $W(g_0,{\varepsilon_1},g_1,{\varepsilon_2}, \ldots, g_{n-1},{\varepsilon_n})$ where  $(g_0,{\varepsilon_1},g_1,{\varepsilon_2}, \ldots, g_{n-1},{\varepsilon_n})$ runs over cyclically reduced sequences with a sequences of $\varepsilon$'s given by  $({\varepsilon_1},{\varepsilon_2}, \ldots,{\varepsilon_n})$.
\end{prop}

 With arguments similar to those of the proof of \propc{filtration mini}, we can prove the following result.
\start{prop}{filtration mini hnn} With the above notations,  the subspace 
$W_{n}(y)$ is the disjoint union of $W(g_0,{\varepsilon_1},g_1,{\varepsilon_2}, \ldots, g_{n-1},{\varepsilon_n})$ where $W(g_0,{\varepsilon_1},g_1,{\varepsilon_2}, \ldots, g_{n-1},{\varepsilon_n})$ runs over all cyclically reduced sequences of $n$ terms of the HNN extension determined by $y$. Therefore, $W_n(y)$ is invariant under $\mathrm{ad}_y$.
\end{prop}

\start{rem}{dehn local hnn}  It is not hard to see that the subspaces  $W(g_0,{\varepsilon_1},g_1,{\varepsilon_2}, \ldots, g_{n-1},{\varepsilon_n})$
 are also invariant 
under the map induced by the Dehn twist around $y$. Indeed, using \lemc{separating representative} , one sees that the Dehn twist around $y$ of the conjugacy class of $g_0,t^{\varepsilon_1},g_1,t^{\varepsilon_2}, \ldots, g_{n-1},t^{\varepsilon_n})$ (where $(g_0,t^{\varepsilon_1},g_1,t^{\varepsilon_2}, \ldots, g_{n-1},t^{\varepsilon_n})$ is a cyclically reduced sequence) is represented by one of the following products:
$$
g_0u_{1}t^{\varepsilon_1}g_1u_{1}t^{\varepsilon_2} \ldots, g_{n-1}u_{n}t^{\varepsilon_n},
$$
where one of the following holds.
\begin{numlist}
\item For each $i \in \{1,2,\cdots,n\}$, if $\varepsilon_{i}=1$ then $u_{i}=y$ and if $\varepsilon_{i}=-1$ then $u_{i}=\varphi(\overline{y})$.
\item For each $i \in \{1,2,\cdots,n\}$, if $\varepsilon_{i}=1$ then $u_{i}=\overline{y}$ and if $\varepsilon_{i}=-1$ then $u_{i}=\varphi({y})$.
\end{numlist}
where the choice  between (1) and (2)  is determined by the orientation of the surface and the orientation of $y$.
\end{rem}

\section{Application: The mapping class group and the curve complex}\label{applications 3}

Let $\Sigma$ be a compact oriented surface.  By $\Sigma_{g,b}$ we denote an oriented surface with genus $g$ 
and $b$ boundary components. If $\Sigma$ is a surface, we denote by $\mathrm{Mod}(\Sigma)$ the \emph
{mapping class group of $\Sigma$}, that is set of homotopy classes of orientation preserving homeomorphisms of $
\Sigma$.  We study automorphisms of the Goldman Lie algebra that are related to the mapping class group. The first two results,Theorems \ref{oriented curves} and \ref{unoriented curves}, apply to all surfaces. The stronger result, \theoc{all oriented}, applies only to surfaces with boundary.

Now we recall the curve complex,  defined by Harvey in \cite{wh}. The \emph{curve complex of $\Sigma$} 
is denoted by $\mathrm{C}(\Sigma)$  is the simplicial complex  whose vertices are  isotopy classes  simple closed 
curves on $\Sigma$ which are not null-homotopic and not homotopic to a boundary component. If $\Sigma \ne 
\Sigma_{0,4} $ and $\Sigma \ne \Sigma_{1,1}$ then a set of $k+1$ vertices of the curve complex is the $0$-
skeleton of a $k$-simplex if the corresponding  minimal intersection number of all pairs of vertices is zero, that is, if 
every pair of vertices have disjoint representatives.

For $\Sigma_{0,4}$ and $\Sigma_{1,1}$ two vertices are connected by an edge when the curves they represent 
have minimal intersection (2 in
the case of $\Sigma_{0,4}$ , and 1 in the case of $\Sigma_{1,1}$).  

The following isomorphism is a theorem of Ivanov
\cite{iv} for the case of genus at least two.  Korkmaz \cite{ko} proved the result for genus at most one and Luo gave 
another proof  that covers all possible genera. \cite{fl}. Our discussion below is based on the formulation of Margalit \cite{margalit}.

\start{theo}{ikl}(Ivanov-Korkmaz-Luo) Let $\Sigma$ be an orientable surface of negative Euler characteristic and let 
$\map{h}{\mathrm{Mod}(\Sigma)}{\mathrm{Aut\hspace{0.05cm} C}(\Sigma)}$ be the natural map. If $\Sigma \ne 
\Sigma_{0,3}$ then
\begin{numlist}\item $h$ is surjective  if and only if $\Sigma \ne \Sigma_{1,2}$ and 
\item  $h$ is injective  if and only if $\Sigma \notin\{ \Sigma_{1,1}, \Sigma_{1,2},  \Sigma_{2,0}, \Sigma_{0,4}\}$.
\end{numlist}
 \end{theo}

\start{theo}{luo}(Luo) Any automorphism of  $C(\Sigma_{1,2})$ preserving the set of vertices represented by 
separating loops is induced by a self-homeomorphism of the surface $\Sigma_{1,2}$. 
\end{theo}

We need the following result from \cite{ck}.
 \start{theo}{ck}(Chas - Krongold) Let $S$ be an oriented surface with non-empty boundary and let $x$ be a free homotopy class of curves in 
$S$. Then 
 $x$ contains a simple representative if and only if   $\langle x, x^3 \rangle =0$, where $x^3$ is the conjugacy class 
that wraps around $x$ three times. \end{theo}

\start{lem}{sep} 
Let $x$  be a free homotopy class of oriented simple closed curves. Then $x$ has  a non-separating representative 
if and only if there exist simple class $y$ such that the minimal intersection number of $x$ and $y$ is equal to one. 
\end{lem}
\begin{proof} If $x$ has a non-separating representative, the existence of $y$ with minimal intersection number 
equal to one can be proved as in  the proof of \lemc{hnn construction}. Conversely, if $x$ contains a separating 
representative, then the minimal intersection number of $x$ and any other class is even.\end{proof}

\start{theo}{oriented curves} Let $\Omega$ be a bijection on the set $\pi^\ast$ of free homotopy  classes of closed 
curves on an oriented surface. Suppose the following
\begin{numlist} 
\item $\Omega$ preserves simple curves.
\item If $\Omega$ is extended linearly to the free $\Z$ module generated by $\pi^\ast$ then $\Omega$ preserves 
the Goldman Lie bracket. In symbols 
$[\Omega(x),\Omega(y)]=\Omega([x,y])$
for all $x, y \in \pi^\ast$.
\item For all $x \in \pi^\ast$, $\Omega(\overline{x})=\overline{\Omega(x)}$.
 \end{numlist}
Then the restriction of  $\Omega$ to the subset of simple closed curves is induced by an element of the mapping 
class group. Moreover, if $\Sigma \notin\{\Sigma_{1,1},\Sigma_{1,2} ,   \Sigma_{2,0}, \Sigma_{0,4}\}$ then the 
restriction of $\Omega$ to the subset of simple closed curves is induced by a unique element of the mapping class 
group.
\end{theo}
\begin{proof} Since $\Omega(\overline{x})=\overline{\Omega(x)}$, $\Omega$ induces a bijective map $\widehat
{\Omega}$ on $\widehat{\pi}=\{x+\overline{x} : x \in \pi^\ast\}$.   

Since $\Omega$ preserves the oriented Goldman bracket and the "change of direction", then it preserves the 
unoriented Goldman bracket. 
Let $x$ be a  class of oriented curves that contains a simple representative, let  $y$ be any class. 

Since $\widehat{\Omega}$ preserves the unoriented Goldman Lie bracket, by \theoc{unoriented},  the minimal 
intersection number of $x$ and $y$ equals the minimal intersection number of $\widehat{\Omega}(x)$ and $
\widehat{\Omega}(y)$. Thus, by \theoc{ikl}, if the surface is not the torus with two holes, $\widehat{\Omega}$ is 
induced by an element of the mapping class group (Observe that the "special" curve complexes are covered by this 
property). Moreover if the surface is not in $\{ \Sigma_{1,1}, \Sigma_{1,2},  \Sigma_{2,0}, \Sigma_{0,4}\}$ then $
\widehat{\Omega}$ is induced by unique element of the mapping class group.

Now we study the case of the torus with two holes. By \lemc{sep},  $x$ is a separating simple closed curve, if and 
only if  $\Omega(x)$ is separating. Thus, $\widehat{\Omega}$ maps bijectively the set of unoriented separating 
simple closed curves onto  itself. Hence, by \theoc{luo},  $\widehat{\Omega}$ is induced  by an element of the 
mapping class group.
\end{proof}

By arguments similar to those of \propc{oriented curves}, one can prove the following.

\start{theo}{unoriented curves} Let $\Gamma$ be a bijection on the set $\widehat{\pi}$ of unoriented free 
homotopy  classes of closed curves on an oriented surface. Suppose the following
\begin{numlist} 
\item $\Gamma$ preserves simple curves.
\item If $\Gamma$ is extended linearly to the free $\Z$ module generated by $\pi^\ast$ then $\Gamma$ preserves 
the unoriented Goldman Lie bracket. In symbols 
$[\Gamma(x),\Gamma(y)]=\Gamma([x,y])$
for all $x, y \in \widehat{\pi}$.
\end{numlist}
Then the restriction of $\Gamma$ to the subset of simple closed curves is induced by an element of the mapping 
class group. Moreover, if $\Sigma \notin\{\Sigma_{1,1},\Sigma_{1,2} ,   \Sigma_{2,0}, \Sigma_{0,4}\}$ then the 
restriction of $\Gamma$ is induced by a unique element of the mapping class group.
\end{theo}

\start{theo}{all oriented} Let $\Omega$ be a bijection on the set $\pi^\ast$ of free homotopy  classes of curves on 
an oriented surface with non-empty boundary. Suppose the following
\begin{numlist} 
\item If $\Omega$ is extended linearly to the free $\Z$ module generated by $\pi^\ast$ then $[\Omega(x),\Omega
(y)]=\Omega([x,y])$
for all $x, y \in \pi^\ast$.
\item For all $x$ in $\pi^{\ast}$, $\Omega(\overline{x})=\overline{\Omega(x)}$ .
\item For all $x$ in $\pi^{\ast}$,  $\Omega({x^3})={\Omega(x)}^3$.
\end{numlist}
Then the restriction of 
$\Omega$ to the set of free homotopy classes with simple representatives is induced by an element of the mapping 
class group. Moreover, if $\Sigma \notin\{ \Sigma_{1,1},   \Sigma_{2,0}, \Sigma_{0,4}\}$ then $\Omega$ is induced 
by a unique element of the mapping class group.
\end{theo}

\begin{proof}
Let $x$ be an oriented closed curve. By hypothesis, $[\Omega(x),\Omega(x)^3]=\Omega([x,x^3])$. Thus, $[x,x^3]
=0$ if and only if $[\Omega(x),\Omega(x)^3]=0$. 
Thus by \theoc{ck}, $\Omega(x)$ is simple if and only if $x$ is simple. Then the result follows from \theoc{oriented 
curves}.
\end{proof}

All these results "support" Ivanov's statement in \cite{iv2}:

\begin{proclama}{Metaconjecture} "Every object naturally associated to a surface $S$ and having a sufficiently rich 
structure has $\mathrm{Mod}(S)$ as its group of automorphisms. Moreover, this can be proved by a reduction 
theorem about the automorphisms of $\mathrm{C}(S)$." 
\end{proclama}

In this sense, the Goldman Lie bracket combined with the power maps, have a "sufficiently rich" structure.

\section{Questions and open problems}\label{questions}

\start{prob}{first rem}
Etingof   \cite{et}  proved using algebraic tools that the center of Goldman Lie algebra of a closed oriented surface consist in the one dimensional subspace generated by the trivial loop. On the other hand, if the surface has non-empty boundary, it is not hard to see that linear combinations of  conjugacy classes of curves parallel to the boundary components are in the center. Hence, it seems reasonable to  conjecture that the center consist on linear combinations of conjugacy classes of boundary components. 
It will be interesting to use the results of this work to give a complete characterization of the center of the Goldman Lie algebra.

If $u$ is an element of the vector space of unoriented curves on a surface, and $u$ is a linear combination of classes that admit a simple representative then $u$ is not in the center of the Goldman Lie algebra. To prove this, we need to combine our results with the fact that given two different conjugacy classes $a$ and $b$ that admit a simple representative, there exists a simple conjugacy class $c$ such that the intersection numbers $i(a,c)$ and $i(b,c)$ are distinct (see \cite{flp} for a proof that such $c$ exists). Nevertheless, by results of Leininger \cite{len} we know that this argument cannot be extended to the cases of elements of the base that only have self-intersecting representatives. 

\end{prob}

\start{prob}{second rem} As mentioned in the Introduction, Abbaspour \cite{hossein} studied whether a three manifold is hyperbolic by means of the generalized Goldman Lie algebra operations. He used free products with amalgamations for this study. We wonder if it would be possible to combine his methods with ours, to give a combinatorial description of the generalized String Topology operations on three manifolds. In this direction, one could study the relation between number of connected components of the output of the Lie algebra operations and intersection numbers. 
\end{prob}

\start{prob}{third rem}   We showed that subspaces $W(w_1, w_2, \ldots, w_n)$  defined in Subsection \ref{double cosets separating} are $\mathrm{ad}_x$-invariant. Let $z$ be a representative of a conjugacy class in  $W(w_1, w_2, \ldots, w_n)$. It would be interesting to define precisely and study the ''number of twists around $x$" of the sequence $(\mathrm{ad}_x^n(z))_{n \in \Z}$ and how this number of twists changes under the action of $\mathrm{ad}_x$. Observe also that the Dehn twist around $x$ increases or decreases the number of twists at a faster rate.  These problems are related to the discrete analog of Kerckhoff's convexity  \cite{Ke} found by Luo \cite{luo2}.
\end{prob}

\start{prob}{fourth rem}  The Goldman Lie bracket of two conjugacy classes, one of them simple, has no cancellation. On the other hand, there are examples (for instance \exc{non-simple}) of pairs of classes with bracket zero and non-zero minimal intersection number. We wonder which is the topological characterization of pairs of intersection points of two curves, for which the corresponding term of the bracket cancel. In other words, what "causes" cancellation. The tools to answer to this question may involve the study Thurston's compactiÞcation of Teichm\"{u}ller space in the context of Bonahon's work on geodesic currents \cite{bonahon}.
\end{prob}

\start{prob}{fifth rem}  Dylan Thurston \cite{DT} proved a suggestive result: Let $m$ be a union of  conjugacy classes of curves on an orientable surface and let $s$ be the conjugacy class of a simple closed curve. Consider representatives $M$ of $m$, $S$ of $s$ which intersect (and self-intersect) in the minimum number of points. Denote by $P$ one of the self-intersection points of $M$. Denote by $M_{1}$ and $M_{2}$ the two possible ways of smoothing the intersection at $P$. Denote by $m_{1}$ and $m_{2}$ respectively the conjugacy classes of $M_{1}$ and $M_{2}$. Then Dylan Thurston's result is
$$
i(m,s) = \max (i(m_{1},s), i(m_{2},s)).
$$
These two ``smoothings'' are the two local operations one makes at each intersection point to find a term of the unoriented Goldman Lie bracket (when the intersection point $P$ is not a self-intersection point of a curve). It might be interesting to explore the connections of Dylan Thurston's result and our work.
\end{prob}

\start{rem}{sequences of double cosets} In a subsequent work we will give a combinatorial description of the set of  cyclic sequences of double cosets under the action of $\mathrm{ad}$. Also, we will study under which assumptions the cyclic sequences of double cosets are a complete invariant of a conjugacy class.
\end{rem}

\bibliographystyle{amsalpha}

\textsc{
Moira Chas,
Department of Mathematics,
SUNY at Stony Brook,
Stony Brook, NY, 11794
}

\emph{E-mail address}{:\;\;}\texttt{moira@math.sunysb.edu}

\end{document}